\documentclass{article}

\usepackage{amsmath}
\usepackage{arxiv}
\usepackage{lscape}
\usepackage[utf8]{inputenc} % allow utf-8 input
\usepackage[T1]{fontenc}    % use 8-bit T1 fonts
\usepackage{hyperref}       % hyperlinks
\usepackage{url}            % simple URL typesetting
\usepackage{booktabs}       % professional-quality tables
\usepackage{amsfonts}       % blackboard math symbols
\usepackage{nicefrac}       % compact symbols for 1/2, etc.
\usepackage{microtype}      % microtypography
\usepackage{lipsum}
\usepackage{graphicx}
\graphicspath{ {./images/} }

\title{Queueing-Inventory Systems: A Survey}

\author{
 Salini K \\
  Department of Mathematics\\
  Central University of Kerala\\
  Kasaragod-671320, India \\
  \texttt{salini.k@cukerala.ac.in} \\
   \And
 Arya  P S \\
  Department of Mathematics\\
  Central University of Kerala\\
  Kasaragod-671320, India \\
  \texttt{arya.p.s@cukerala.ac.in} \\
  \And
 Manikandan Rangaswamy* \\
  Department of Mathematics\\
  Central University of Kerala\\
  Kasaragod-671320, India  \\
  \texttt{mani@cukerala.ac.in} \\
}

\begin{document}
\maketitle
\begin{abstract}
Queueing-inventory systems are integrated systems consisting of two emerging fields in applied probability, namely ``Queues" and ``Inventory". In this paper, we present a comprehensive review of the theory and applications of queueing-inventory systems from 2016 to 2023 and a detailed discussion of the recent literature by encompassing the topics like product-form solutions, production-inventory, perishable inventory, the retrial of unsatisfied customers/vacation/vacation interruption, multi-servers, multi-commodity, multi-class customers, postponed demands, and game-theoretic analysis on queueing-inventory systems, and other applications in various fields such as energy harvesting in the wireless network, supply chain, etc. The review provides insights into the importance of queueing-inventory systems, systematically discussing the developments and their applications, challenges, and future directions.
\end{abstract}

% keywords can be removed

\section{Introduction}
Mathematical inventory theory had its origin in the development of the EOQ (Economic Order Quantity) formula by Harris \cite{har} in 1913. Harris’s work was under a deterministic setup. Subsequently, several researchers started investigating problems in deterministic inventory. The impetus for stochastic models in inventory is due to Arrow \textit{et. al} \cite{arrow} in 1958. This leads to the analogue of the EOQ formula in stochastic inventory models. Moran \cite{pap} and Prabhu (\cite{pnu1}, \cite{pnu2}) were the pioneers in this area of research. For basic deterministic and Stochastic models in inventory, one may refer to Hadley and Whitin \cite{hw} and Naddor \cite{naddor}. All work reported in inventory theory until the early 1990s made the very strong assumption that the time required to serve customers was negligible. In most situations, this assumption is unrealistic. The first work reported on the inventory with positive service time, also called the ``queueing-inventory system," was by Melikov and Molchanov \cite{mam}, Sigman and Simchi-Levi  \cite{ssl} in 1992. 

In contrast to a traditional queuing system, the time for serving an item is considered negligible. Still, most real-life situations need more time to serve the items in the inventory. Hence, at first glance, an inventory with a positive service time may seem to be a classical queueing problem. The classical queueing model, however, assumes that if customers are available and the server is ready to serve, service could be provided. In other words, under the classical queueing theory, the quantity of the items needed for service is considered or ignored. But, in reality, serving an item to a customer from the inventory will take a random amount of time. Furthermore, the items in the inventory must be available for service in a queueing-inventory system in addition to the customers and servers. However, there's always a chance that the inventory won't be available in stock. The customer leaves with one or more items at the service completion epoch. As a result, the inventory level decreases, corresponding to the number of customers the system serves. Also, it is obvious that in a situation of positive lead time, the server may be idle even if the customers are waiting for items in the queue when the items are not in stock.

 Researchers have closely studied queueing systems connected with inventories for the last three decades. After the pioneering work by  Melikov and Molchanov \cite{mam}, Sigman and Simchi-Levi  \cite{ssl}, several articles appeared in various reputed journals concerning with QISs by researchers. An extensive literature review on queueing-inventory systems can refer to Krishnamoorthy \textit{et al.} \cite{klm}, and \cite{ksn} for a quick, albeit incomplete, work up to 2018. In addition, a detailed survey of the queueing-inventory models with phase-type service distributions from 1985 to 2016 is provided by Kyung Hwan Choi and Bong Kyoo Yoon \cite{ckh}. Apart from those review articles, this paper encompassed the articles on queueing-inventory systems from 2017-2023, and also we included the articles which are not listed in the previous surveys. The product-form solutions will be considered an achievement in any stochastic system. In this survey, we narrate the research developments/contributions on queueing-inventory systems with different scenarios, including the articles on the game-theoretic queueing-inventory systems. 
 
 This review article is organized as follows: section 2.1 considers a queueing-inventory system with a product-form solution and provides a brief discussion on each model. Sections 2.2, 2.3, and 2.4 provide an overview of the work related to single-server, multi-server queueing-inventory systems, and queueing-inventory systems with heterogeneous customers, respectively. Section 2.5 examines multi-commodity queueing-inventory systems. In sections 2.6, 2.7, and 2.8, we investigated queueing-inventory systems with perishable inventory, postponed demands, and retrials. Section 2.9 summarizes the work done on the production-inventory system. The queueing-inventory system with the random environment is discussed in section 2.10. The works on game-theoretic queueing-inventory systems are included in section 2.11. The article's on queueing-inventory systems with different contexts are given in section 2.12. The conclusion and future work are discussed in section 3. Tables 1- 7 provide an overview of the work reported in queueing-inventory systems from 2016-2023, and the classification of the Game Theoretic queueing-inventory system is shown in Table 8. The tables are classified based on keywords.

 \subsection*{\textbf{Abbreviation used in the sequel are:}}
	\begin{itemize}
		\item QMCD: Queueing and Markov Chain Decomposition
		\item MAM: Matrix-Analytic Method
		\item QIS: Queueing-Inventory System
            \item PIS: Production-Inventory Systems
		\item MAP: Markovian Arrival Process
		\item MMAP: Marked Markovian Arrival Process
		\item BMAP: Batch Markovian Arrival Process
		\item PH: Phase Type
            \item RE: Random Environment
            \item CLT: Common LifeTime
		\item SMA: State Merging Algorithm
            \item SCV: Squared Coefficient of Variation
            \item MTS: Make-to-Stock
            \item MTO: Make-to-Order
	\end{itemize}
 
 The terminology and terms of queueing-inventory systems are adopted from Krishnamoorthy \textit{et al.} \cite{ksn}. Also, a few definitions/concepts/phrases that are essential for the game-theoretic approach of the queueing-inventory system are provided, and it may help the readers better understand the game-theoretic perspective of queueing-inventory models.    \\ \\

  \begin{itemize}
      \item \textit{Lead Time}: The lead time is the time lag between the epoch at which an order for inventory replenishment is placed and the epoch at which the order is materialized. \\
  \item  \textit{$(s, S)$-policy:} This is an inventory control policy in which the replenishment quantity is equal to the amount required to return the level to $S$ at the replenishment epoch. The reorder level is $s$, and the maximum number of items stored in $S$.  \\
   \item \textit{$(S-1, S)$-policy:} This control policy requires placing an order whenever an item is sold out. It is often referred to as the one-for-one ordering policy. This policy is recommended for slow-moving and expensive items. \\
  \item \textit{$(s, Q)$-policy:} The number of items ordered in this control policy is fixed and equals $Q = S - s$. Here, s represents the reorder level, and $Q$ represents the fixed order quantity. To avoid perpetual order placement for replenishment, $Q$ is made larger than $s$.  \\
 \item In \textit{random order size}, the decision of the order size is according to a discrete probability function p on the set $\{1, 2, ..., S\}$. The size of a replenishment order is $k$ with probability $p_{k}$, such that $\sum_{k=1}^{S} p_{k} = 1$.   \\
 
  \item \textit{Phase-type distribution} (Neuts \cite{nmf}): The continuous phase-type distributions are introduced as a natural generalization of the exponential and Erlang distributions. A phase-type distribution is obtained as the distribution of the time until absorption in a Markov chain with a finite state space and an absorbing state.\\ \\
  Consider a Markov process $\varphi=\{X(t):t\geq 0\}$ with a finite state space $\{1,2,\dots,m+1\}$ and the infinitesimal generator matrix $Q=\begin{pmatrix}
     T & T^{0} \\ 
     0 & 0 
    \end{pmatrix}$ 
    where $T$ is a square matrix of order $m$, $T^{0}$, a column vector and $0$ the zero row vector of the same dimension. The elements of the matrices $T$ and $T^{0}$ satisfy $T_{ii} < 0$ for $1 \leq i \leq m$, $T_{ij} \geq 0$ for $i\neq j$; $T_{i}^{0} \geq 0$ and $T_{i}^{0} > 0$ for at least one $i$, $1 \leq i \leq m$, $T_{0}$ $i$ is the $i^{th}$ element of $T_{0}$and $T_{e} + T^{0} = 0$. \\ \\
    The initial distribution of $\varphi$ is given by the row vector $(\alpha, \alpha_{m+1})$ where $\alpha$ a row vector of dimension $m$ with the property that $\alpha e+\alpha_{m+1} = 1$. The states $\{1, 2, ...,m\}$ shall be transient, while the state $m + 1$ is absorbing. \\ \\ 
    Let $Z:= inf\{t \geq 0: X(t) = m + 1\}$ denote the random variable of the time until absorption in state $m+1$. The distribution of $Z$ is called phase-type distribution (or shortly PH distribution) with representation $(\alpha, T)$. The dimension $m$ is called the order of the distribution. The states $\{1, ...,m\}$ are also called phases, which gives rise to the name phase-type distribution. The distribution function of $Z$ is given by 
   \begin{equation*}
       F(t):=P(Z\leq t)=1-\alpha \textbf{e}^{Tt}e, \quad for\ all \ t\geq 0 
   \end{equation*} and the density function is 
   \begin{equation*}
       f(t) = \alpha \textbf{e}^{Tt} T^{0}, \quad t>0.
   \end{equation*}   
 \\
\item \textit{Markovian arrival process:} Markovian arrival process $(MAP)$ introduced in Neuts \cite{nmf}is a rich class of point processes that includes many well-known processes such as Poisson, PH-renewal processes, and Markov-modulated Poisson process. A significant feature of the $MAP$ is the underlying Markovian structure that fits ideally in the context of matrix-analytic solutions to stochastic models. $MAP$ is a generalization of the Poisson process, which keeps many useful properties of the Poisson process. For example, the memoryless property of the Poisson process is partially preserved by the $MAP$ by conditioning on the phase of the underlying Markov chain. Any stochastic counting process can be approximated arbitrarily closely by a sequence of Markovian arrival processes. $MAP$ is a convenient tool to model both renewal and non-renewal arrivals. Chakravarthy \cite{chy} provides an extensive survey of the batch Markovian arrival process $(BMAP)$ in which arrivals are in batches whereas it is in singles in $MAP$.\\ \\
    A continuous-time Markov chain $\{(N(t), I(t)), t \geq 0\}$ with state space \\ ${(i, j) : i \geq 0, 1 \leq j \leq m} $ and infinitesimal generator
    $Q= \begin{pmatrix}
     D_{0} & D_{1} \\ 
     &    D_{0} & D_{1}\\
      &  & . &. \\
      & & & . &.\\
      & & & & .&.
    \end{pmatrix}$ is called an $MAP$ with matrix representation $(D_{0},D_{1})$. Here,$D_{0}$ and $D_{1}$ are square matrices of order $m$ where $m$ is a positive integer. The diagonal elements of $D_{0}$ are negative and its off-diagonal elements are non-negative, $D_{1}$ has all its elements non-negative, and $D_{0} + D_{1}$ is an infinitesimal generator. Let $D_{0}=\left( d_{ij}^{(0)}\right)$ and $D_{1}=\left( d_{ij}^{(1)}\right)$, then $d_{ij}^{(0)}$ is the rate of transitions from phase $i$ to $j$ without an arrival, for $i\neq j$; $d_{ij}^{(1)}$ is the rate of transitions from phase $i$ to $j$ with an arrival and $-d_{ii}^{(0)}$ is the total rate of events in phase $i$. Let $N(t)$ denote the number of arrivals in $(0, t)$ and $I(t)$ the phase of the Markov chain at time $t$. Let $\pi$ be the stationary probability vector of $D$. Then the constant $\lambda = \pi D_{1}e$, referred to as \textit{fundamental rate}, gives the expected number of arrivals per unit of time in the stationary version of the $MAP$.
    \\
    \item   \textit{Blocking set:} If it is assumed that the system will not admit new arrivals as long as it is in a specified subset of the process's state space, that subset is referred to as a blocking set for new arrivals. Arrivals that occur while the system is in the blocked set are considered to be lost by the system. \\
  \item  \textit{Perishability:} The term perishability refers to the ageing of stored items. Things perish one at a time. \\
 
   \item \textit{Strategy:} a strategy is a set of choices or actions a player takes to achieve a certain outcome in a game. A player's strategy depends on their understanding of the game, the available information, and the behavior of other players. A strategy can be either pure or mixed. A \textit{pure strategy} is a specific choice of action that a player makes with certainty. In contrast, a \textit{mixed strategy} is a probability distribution over a player's actions.\\
   \item \textit{Equilibrium strategy:} an equilibrium strategy is a set of strategies that result in a stable outcome when all players in the game follow those strategies. The most common type of equilibrium strategy in game theory is the \textit{Nash equilibrium}, which is named after mathematician \textit{John Nash}. A Nash equilibrium is a set of strategies where no player can improve their outcome by unilaterally changing their strategy, assuming that all other players are following their equilibrium strategies.\\
   \item \textit{Observable queueing-inventory system:} In an observable queueing-inventory system, customers can be directly observed or measure the queue length and the number of items available in the inventory by the system or the service provider.\\
    \item \textit{Unobservable queueing-inventory system:} In an unobservable queueing-inventory system, customers cannot be directly observed or measure the queue length and the number of items available in the inventory neither by the service provider nor by the system. \\ 
   \item \textit{Duopoly:} Duopoly is a term used in game theory to describe a market structure where there are only two firms that produce and sell a particular product. In a duopoly, the actions of one firm can have a significant impact on the profits of the other firm, creating a complex strategic interaction between the two firms.\\ 
   \item \textit{Payoff:} a payoff is a numerical value that represents the reward or benefit that a player receives as a result of the outcome of a game.\\ 
   \item \textit{Make-to-Order:} Make-to-order is a production strategy where products are only manufactured once a customer places an order. In this method, the manufacturing process begins after the order is received, and the product is tailored to meet the customer's specific requirements.\\ 
   \item \textit{Make-to-Stock:} Make-to-stock is a manufacturing strategy where products are produced in advance of customer demand based on sales forecasts or historical data.
   \end{itemize}

\section{Queueing-Inventory Models}
\label{sec: Queueing-Inventory Models}
Customers arrive at the service facility one by one and require service in a queueing-inventory system. An item from the inventory is needed to complete the customer service. A served customer departs immediately from the system, and the on-hand inventory decreases by one at the moment of service completion. An outside supplier provides the inventory. According to Schwarz \cite{scz} \textit{et al.}, this system is known as a queueing-inventory system.

\subsection{\textbf{Queueing-Inventory system with a product-form solution}}

Baek \textit{et al.} \cite{jwbaek} is the first reported work on a queuing system with a continuous-type $(s, Q)$ control policy, lost sales and general lead time. The authors  considered an $M/M/1$ queuing model with an attached continuous-type inventory system motivated by chemical and 3D printing. They derived stationary joint probability distribution of queue length and inventory level explicitly. In this model, at the service completion time, each customer requires an exponentially distributed random amount of inventory ($H$); the inventory level decreases by the stock requirement  if the on-hand inventory is sufficient. The customers arrive according to the Poisson process at rate $\lambda$, and $\mu$ is the independent, identically and exponentially distributed service times. The stationary joint probability distribution of the corresponding model is,
	\begin{equation*}
		{\Tilde{\psi}}^D_n (x,w) = P_{M/M/1}(n).{\Tilde{X}}^D(x,w), \quad (n=0,1,2,\dots, w\geq 0, 0\leq x <s),
	\end{equation*}
	\begin{equation*}
		\psi^R_n (x) = P_{M/M/1}(n).{X}^R(x), \quad (n=0,1,2,\dots, s\leq x\leq s+Q),
	\end{equation*}
	\begin{equation*}
		A^R_n = P_{M/M/1}(n).B^R, \quad (n=0,1,2 \dots),
	\end{equation*}
	where \begin{equation*} 
		P_{M/M/1}(n) = \left(1-\frac{\lambda}{\mu}\right) \left(\frac{\lambda}{\mu}\right)^{n}, \quad (n=0, 1, \dots),
	\end{equation*}
	\begin{eqnarray*}\small
		\Tilde{X}^{D}(x,w)= \mathop{lim}\limits_{t\rightarrow \infty} Pr[I(t)\leq x, w<V_{+}(t)\leq w+dw, \xi(t)=1], \\ (w\geq 0, 0 \leq x<s),
	\end{eqnarray*}
	\begin{equation*}
		\Tilde{X}^{R}(x)= \mathop{lim}\limits_{t\rightarrow \infty} Pr[I(t)\leq x, \xi(t)=0], \quad (s \leq x< s+Q),
	\end{equation*}
	\begin{equation*}
		\Tilde{B}^{R} = \mathop{lim}\limits_{t\rightarrow \infty} Pr[I(t)=Q,\xi(t)=0].
	\end{equation*}
	Finally, the authors achieved to decompose the steady-state probability distribution as the stationary joint probability distribution of their proposed model and modified model.

	An $M/M/1$ QIS with batch demands and lost sales considered by Yue D \textit{et al.} \cite{yzq} discussed the two cases of lost sales: the partial-lost sales and the full-lost sales models. In the partial lost sales model, the customer takes away all the items in the inventory. The unsatisfied demands are lost if the on-hand inventory is less than the size of the batch demands of an arrived customer. While in the second model, the customer leaves without taking any item from the inventory, and all customer demands are lost completely if the on-hand inventory is less than the size of the batch demands of an arrived customer. Using the QBD process, they obtained the stationary distributions of the joint queue length product-form and the on-hand inventory process for both models. The authors derived the steady-state probability distribution for the first model with the stability condition, $ \rho = \frac{\lambda}{\mu} < 1$ as follows:
		\begin{equation*}
		\pi_{n} = (1-\rho) \rho^n x_i, \quad n > 0,
	\end{equation*}
	where the first part of the above, $(1-\rho) \rho^n$, $(n \geq 0) $ is the joint stationary probability distribution of the queue length process as a typical $M/M/1$ queueing system with the parameters $\lambda$ and $\mu$, and the second part, $x_i$ $(0 \leq i \leq S)$
	is the joint stationary probability distribution of the on-hand inventory process. 
	Similarly, the stationary probability distribution for the second model is derived as follows,
	\begin{equation*}
		\pi_{n} = (1-\rho) \rho^n \mbox{\boldmath {{$y$}}}, \quad n \geq 0,
	\end{equation*}  where \mbox{\boldmath {{$y$}}}$= (y_{0}, y_{1},\dots, y_{S})$ be the stationary distribution of the process  \\ $\left\{I(t), t\geq 0\right\} $, $\lambda$ is the arrival rate of customers follows Compound Poisson process. The service time is exponentially distributed with the parameter $\mu$. They proved that the stationary distribution of the joint queue length and the on-hand inventory process for each model has a product-form by assuming arriving customers are lost when the inventory level is zero; they might not obtain the product-form.

	PIS models by considering both a positive service time and vacations of a production facility is a new approach discussed by Yue D \textit{et al.} \cite{yq}. The system has a single production facility that produces homogeneous items with an $(s, S)$ control policy. The production facility takes a vacation of random duration once the inventory level reaches the level $S$. They obtained  the product-form solution of the stationary joint distribution of the queue length and the on-hand inventory level by using the direct method, and it is of the form: 
	\begin{equation*}
		P_i = \left(1- \frac{\lambda}{\mu}\right) \left(\frac{\lambda}{\mu}\right)^n \mbox{\boldmath {{$\pi$}}}, \quad i \geq 0,
	\end{equation*} where $\lambda$ is the arrival rate of customers follows the Poisson process, the service times are exponentially distributed with mean $\dfrac{1}{\mu}$ and the steady-state probability vector, \mbox{\boldmath {{$\pi$}}}$ = (\pi(0,0), \pi(0,1), \pi(1,0), \pi(1,1), \dots, \pi(S-1,0), \pi(S-1,1), \pi(S,0))$ is: 
	\begin{equation*}
	    \pi(0,0) = \frac{\lambda}{\theta} \left(\frac{\lambda}{\lambda+\theta}\right)^{s} \pi(S,0),
	\end{equation*}
	\begin{equation*}
	    \pi(i,0) =  \left(\frac{\lambda}{\lambda+\theta}\right)^{s-i+1} \pi(S,0), \ 1\leq i \leq s;
	\end{equation*}
	\begin{equation*}
	    \pi(i,0) = \pi(S,0), \ s+1 \leq i \leq S-1,
	\end{equation*}
	\begin{eqnarray*}\small
        \pi(i,1) = \Bigg\{\frac{\lambda}{\eta-\lambda} \left[1-\left(\frac{\lambda}{\eta}\right)^{S-s}\right]\left(\frac{\lambda}{\eta}\right)^{s-i} +\frac{\lambda}{\lambda+\theta-\eta}\left(\frac{\lambda}{\eta}\right)^{s-i}  -  \frac{\lambda}{\lambda+\theta-\eta}\left(\frac{\lambda}{\lambda+\theta}\right)^{s-i}\Bigg\} \pi(S,0), \\ 0\leq i\leq s;
	\end{eqnarray*}
	\begin{equation*}
	    \pi(i,1) = \frac{\lambda}{\eta-\lambda} \left[1-\left(\frac{\lambda}{\eta}\right)^{S-i}\right] \pi(S,0),  \ s+1 \leq i \leq S-1,
	\end{equation*} where,  \small
        
          \begin{eqnarray*} 
	  \pi(S,0) =  \Bigg\{ \left(\frac{\lambda}{\lambda-\eta}\right)^{2} \left(\frac{\lambda}{\eta}\right)^{S} - \frac{\lambda^{2}\theta}{(\lambda-\eta)^{2}(\lambda+\theta-\eta)}\left(\frac{\lambda}{\eta}\right)^{s} + \frac{\lambda^{2}}{\theta(\lambda+\theta-\eta)} \left(\frac{\lambda}{\lambda+\theta}\right)^{s} - \frac{\eta}{\lambda-\eta} \left(S-s+\frac{\lambda}{\theta}\right)\Bigg\}^{-1}.
	\end{eqnarray*}  \normalsize

	They also derived the performance measures explicitly, constructed a cost function based on some performance measures, and discussed the monotonic behavior of the optimal policy and the optimal cost on the system parameters.

	Yue D \textit{et al.} \cite{ywz} considers a PIS with a service facility, production interruptions, and an $(s, S)$ production policy. The production is interrupted for a random vacation time when the inventory level becomes $S$, with exponentially distributed vacation time of rate $ \theta$. Customer arrivals and service follow the Poisson process and exponential distribution with rate $\lambda$ and $\mu$, respectively. The production is switched on immediately and continues in the mode until the inventory level becomes $S$ If the inventory level depletes to $s$ after vacation. Otherwise, the production facility takes another vacation of random duration (multiple vacations) at the same distribution rate as before when it returns from the vacation if the inventory level exceeds $s$. Using the matrix analytic method, they obtained the product-form solution for the stationary joint distribution of the number of customers and the on-hand inventory level, ie, if $\rho < 1$ ($\rho=\frac{\lambda}{\mu}$)
	\begin{equation*}
		P_{i} = (1-\rho)(\rho)^{i} \mbox{\boldmath{{$\pi$}}} , \quad i=0, 1, \dots
	\end{equation*}
	where $\mbox{\boldmath{{$\pi$}}} = (\pi(0, 0), \pi(0,1), \pi(1, 0), \pi(1, 1), \dots, \pi(S-1, 0), \pi(S-1, 1), \pi(S, 0))$ that is, 
	\begin{equation*}
	    \pi(i,0) = \alpha_{i} \pi(S,0), \quad i= 0, 1, \dots, S-1,
	\end{equation*}
	\begin{equation*}
	    \pi(i,1) = \beta_{i} \pi(S,0), \quad i= 0, 1, \dots, S-1,
	\end{equation*}
	and 
	\begin{equation*}
	    \pi(S,0) = \Bigg[1+ \sum_{i=0}^{S-1}(\alpha_{i}+\beta_{i})\Bigg]^{-1}
	\end{equation*}
	where 
	\begin{equation*}
	    \alpha_{0} = \frac{(\lambda+\gamma)}{\theta}\alpha_{i},
	\end{equation*}
	\begin{equation*}
	    \alpha_{i} = \alpha_{s+1} \prod_{k=i}^{s} \frac{\lambda+(k+1)\gamma}{\lambda+k\gamma+\theta}, \quad i= 1,2,\dots,s,
	\end{equation*}
	\begin{equation*}
	    \alpha_{i} = \alpha_{s+1} \prod_{k=i}^{S-1} \frac{\lambda+(k+1)\gamma}{\lambda+k\gamma}, \quad i= s+1, s+2,\dots, S-1
	\end{equation*}
	and for $i=0,1,\dots, s-1,$ $\beta_{i}$ is determined by
	\begin{equation*}
	    \beta_{i} = \frac{\lambda+(i+1)\gamma}{\eta}\beta_{i+1} + \frac{\lambda+S\gamma}{\eta}-\frac{\theta}{\eta}\sum_{k=i+1}^{s}\alpha_{k},
	\end{equation*}
	\begin{equation*}
	   \beta_{i} = \frac{\lambda+S\gamma}{\eta}\sum_{j=1}^{S-i-1} \left\{1+\frac{1}{\eta^{k}}\prod_{k=1}^{j}[\lambda+(i+k)\gamma]\right\}, \quad i= s, s+1, \dots, S-2,
	\end{equation*}
	\begin{equation*}
	    \beta_{S-1} = \frac{\lambda+S\gamma}{\eta}
	\end{equation*}
	 and $\mbox{\boldmath{{$P$}}} = (P_0, P_1, \dots)$ is the steady-state probability vector.
	They also obtained the cost function explicitly and numerically illustrated some optimization methods, such as a genetic algorithm to obtain the minimal policy $(s^{*}, S^{*})$ to minimize the cost function.

	Krishnamoorthy \textit{et al.}  \cite{kavr} considered an $(s, S)$ production-inventory model with lead time involving local purchases to ensure customer satisfaction and goodwill. The problem  was modeled as a Continuous Time Markov Chain and obtained stochastic decomposition of system states. Analysis of this model has high social relevance as it helps reduce the total expected cost of the production process, which improves the profit in the cottage industry. They obtained steady-state probability vector explicitly as follows:
	\begin{equation*}
		z^{(i)} = (1-\rho) \rho^{i} \mbox{\boldmath{{$\Delta$}}}, \quad i\geq 0,
	\end{equation*}
	where $\bar{Z}$ is the steady-state probability vector of the generator $\hat{A}$ and $\bar{Z}= (z^{0}, z^{1}, z^{2}, \dots)$ and the sub-vectors $z^{i}= (z^{(i,1)}, z^{(i,2)}, z^{(i,3)}, \dots, z^{(i,s)}, z^{(i,s+1,0)}, z^{(i,s+1,1)}, \dots, z^{(i,S-1,0)}, z^{(i,S-1,1)}, z^{(i,S)}), \quad i=0, 1, 2, \dots$ and  $\mbox{\boldmath{{$\Delta$}}} = (\pi_{1}, \pi_{2}, \dots, \pi_{s-1}, \pi_{s}, \Tilde{\pi_{s+1}}, \Tilde{\pi_{s+2}}, \dots, \Tilde{\pi_{S-1}}, \Tilde{\pi_{S}})$ is the steady-state probability vector with negligible service time and $\rho = \frac{\lambda}{\mu}$, $\lambda$ is the rate of arrival of demands follows Poisson process and exponentially distributed processing time with parameter $\mu$. The convexity of the cost function was found numerically and obtained several performance measures.

	Shajin \textit{et al.} \cite{skm} analyzed a single server QIS with infinite capacity, and partial and complete blocking sets of states follow the $(s_1, Q)$ order policy. New arrivals are restricted when the inventory level enters a partially blocking set, while in a complete blocking set other than replenishment, block all activities. Stochastic decomposition is not possible in a partial blocking set, and a product-form solution is obtained only in a complete blocking set, and this is of the form,
	\begin{equation*}
		\eta_i = \left(1-\frac{\lambda}{\mu}\right) \left(\frac{\lambda}{\mu}\right)^{i} \mbox{\boldmath{{$\xi^{'}$}}},  \quad i \geq 0, 
  \end{equation*} 
  ~\\ under the condition $\lambda<\mu$, where $\lambda$ is the arrival rate following the Poisson process, $\mu$ is the exponentially distributed service rate and $ \mbox{\boldmath{{$\xi^{'}$}}}=(\xi_{0}^{'}, \xi_{1}^{'}, \dots \xi_{S}^{'} )$ is the steady-state vector of the infinitesimal generator $W$,
	\begin{equation*}
		\xi_{i}^{'} = \begin{cases}
			0 & 0\leq i\leq s_{2}-1\\
			\frac{\beta}{\lambda}(\frac{\beta+\lambda}{\lambda})^{i-(s_{2}+1)} \xi_{s_{2}}^{'} & s_{2}+1\leq i\leq s_{1}, \\
			\frac{\beta}{\lambda}(\frac{\beta+\lambda}{\lambda})^{s_{1}-s_{2}} \xi_{s_{2}}^{'} & s_{1}+1\leq i\leq Q+s_{2}, \\
			\frac{\beta}{\lambda}(\frac{\beta+\lambda}{\lambda})^{i-(Q+s_{2}+1)} [(\frac{\beta+\lambda}{\lambda})^{(Q+s_{1}+1)-i} -1] \xi_{s_{2}}^{'} & Q+s_{2}+1\leq i\leq Q+s_{1},
		\end{cases}
	\end{equation*}
	where \begin{equation*}
		\xi_{s_{2}^{'}}=    \left[1+Q\frac{\beta}{\lambda}\left(\frac{\beta+\lambda}{\lambda}\right)^{s_{1}-s_{2}}\right]^{-1}.
	\end{equation*}
 The product-form solution is possible even if the lead time follows the general distribution. Numerical illustrations give support for the impact of various performance measures.

	Otten \textit{et al.} \cite{osk} investigated a class of separable systems consisting of parallel production systems at several locations associated with local inventories under a base-stock policy connected with a supplier network. The production system manufactures according to customers' demands on a make-to-order basis. The authors studied two types of lost sales based on local inventory and available inventory. They obtained the product-form stationary distribution where the system assumed to be ergodic, i.e., $\mathop{lim}\limits_{t\rightarrow \infty} P(Z(t)= (n, k))= \pi(n, k)= \xi(n) \cdot \theta(k)$. The customers arrive one by one at location $j$ according to a Poisson process with rate $\lambda_{j}$ and $\mu_{j}$ is the service intensity,  if there are  $l> 0$ orders present, the service intensity is $\nu_{m}(l)$, the product structure of the stationary distribution
	\begin{equation*}
		\xi(n)=\prod_{j\in \bar{J}} \xi_{j}(n_j)=\prod_{j\in \bar{J}} C_j^{-1} \prod_{l=1}^{n_j}\left(\frac{\lambda_j}{\mu_j(l)}\right), \quad n_{j} \in N_{0};
	\end{equation*} and 
	\begin{equation*}
		\theta(k) = C_{\theta}^{-1} \prod_{j\in\bar{J}}\left(\frac{1}{\lambda_{j}}\right)^{\#k_{j}} \cdot \prod_{m\in\bar{M}}\prod_{l=1}^{\#k_{m}}\left( \frac{1}{\nu_{m}(l)}\right), \quad \textbf{k} \in K;
	\end{equation*} and normalization constants,
	\begin{equation*}
		C_{j} = \sum_{n_{j}\in N_{0}} \prod_{l=1}^{n_{j}} \left(\frac{\lambda_j}{\mu_{j}(l)}\right),\ and \ 
		C_{\theta} = \sum_{k\in K} \prod_{j\in\bar{J}}\left(\frac{1}{\lambda_{j}}\right)^{\#k_{j}} \cdot \prod_{m\in\bar{M}}\prod_{l=1}^{\#k_{m}}\left( \frac{1}{\nu_{m}(l)}\right),
	\end{equation*}
	is frequently discovered for typical Jackson networks and relatives. Also, the authors have discussed a significant distinction between the "queueing network" and the regular Jackson networks. As long as users are present at each node, servers in Jackson networks remain constantly busy. However, servers in this production network will be idle while the customers are waiting for service until replenishment occurs. The production system manufactures according to customers' demands on a make-to-order basis. Cost analysis is a separable optimization problem due to decoupling the marginal distributions in a stationary state. Finally, they obtained the explicit product-form  solutions for both stationary distributions, which enabled them to compare the main performance metrics of the systems under the different lost sales admission regimes.

	An $M/M/1$ retrial queue with a storage system under $(s, S)$ policy studied by Shajin and Krishnamoorthy \cite{sdk}. The authors assume that the external arrivals enter directly into an orbit when the server is idle, and the interval between two subsequent retrials follows an exponential distribution. The customer at the head of the orbit is only permitted to access the server. In addition, they assume a strong condition that there is no customer could join in the orbit when the storage system is exhausted or the server is busy. They have constructed a partially blocking set for the arrival process and also provided a geometric distribution, except for the multiplicative constants for the modified retrial queueing process introduced to obtain the stochastic decomposition on the steady-state probability distribution. Since it restricted the entry of the primary customers to the orbit \textit{i.e.}, they are not allowed to join the orbit when a service is proceeding. Only when the server is idle the primary customers are allowed to enter the orbit, and it forms a queue. Another crucial assumption they made is that when the inventory level is zero, no primary customer joins the orbit, or the orbital customers retry; even if they do so, they return to orbit. With all these assumptions in hand, for $\lambda<\theta$ where $\theta$ is the exponentially distributed rate of an interval between two successive repeated attempts, they obtained the steady-state probability distribution as follows: \begin{equation*}
		p_{n}(k,i) =  \begin{cases}
			\left(1-\frac{\lambda}{\theta}\right) \left(\frac{\lambda}{\theta}\right)^{n}\psi(k,i) & \text{for }k=0, 0\leq i\leq S,\\    
			\left(1-\frac{\lambda}{\theta}\right) \left(\frac{\lambda}{\theta}\right)^{n}\psi(k,i) & \text{for }k=1, 1\leq i\leq S;    
		\end{cases}
	\end{equation*} 
	where 
	\begin{equation*}
		\psi(k, i) =\begin{cases}
			\frac{\gamma}{\lambda} \left(\frac{\gamma+\mu}{\mu}\right)^{i} \left(\frac{\gamma+\lambda}{\lambda}\right)^{i-1}
			\psi(0,0), & k=0, 1\leq i \leq s, \\
			\frac{\gamma}{\lambda} \left(\frac{\gamma+\mu}{\mu}\right)^{s} \left(\frac{\gamma+\lambda}{\lambda}\right)^{s}
			\psi(0,0), & k=0, s+1\leq i \leq S-1, \\
			\left[\frac{\gamma}{\lambda}\left(\frac{\gamma+\mu}
			{\mu}\right)^{s} \left(\frac{\gamma+\lambda}{\lambda}\right)^{s} + \frac{\gamma}{\gamma+\lambda+\mu}
			\left[1-\left(\frac{\gamma+\mu}{\mu}\right)^{s} \left(\frac{\gamma+\lambda}{\lambda}\right)^{s} \right]\right] \psi(0,0), & k=0, i=S, \\
			\frac{\gamma}{\mu} \left(\frac{\gamma+\mu}{\mu}\right)^{i-1} \left(\frac{\gamma+\lambda}{\lambda}\right)^{i-1}
			\psi(0,0), & k=1, 1\leq i \leq s,\\
			\frac{\gamma}{\mu} \left(\frac{\gamma+\mu}{\mu}\right)^{s} \left(\frac{\gamma+\lambda}{\lambda}\right)^{s}\psi(0,0), & k=1, s+1\leq i \leq S,
		\end{cases}
	\end{equation*} 
	with
	\begin{equation*}
		\psi(0, 0) = \left\{\left(\frac{\gamma+\mu}{\mu}\right)^{s} \left(\frac{\gamma+\lambda}{\lambda}\right)^{s} (\lambda+\mu)\left[\frac{1}{\gamma+\lambda+\mu} + (S-s)\frac{\gamma}{\lambda\mu}\right] + \frac{\gamma}{\gamma+\lambda+\mu}\right\}^{-1}.
	\end{equation*}
	Under the necessary and sufficient stability condition $\lambda<\theta$, they obtained the steady-state probability vector of the CTMC as 
	$$
		p_{n} = \left(1-\frac{\lambda}{\theta}\right)\left(\frac{\lambda}{\theta}\right)^{n} \mbox{\boldmath{{$\psi$}}},  \quad \text{for} \quad n\geq 0.
  $$
	They derived the stationary joint distribution of the queue length and the on-hand inventory in an explicit product-form solution. They ignored the inventory status by considering that the materials were abundantly available and examined only the status of the server and the number of customers in orbit. For the retrial queue with the steady-state probability vector, the system state equation is,
	\begin{equation*}
		P(i \ customers\ in\ the\ orbit\ and\ server\ is\ busy) = \left( 1-\frac{\lambda}{\theta}\right) \left( \frac{\lambda}{\theta}\right)^{i} \left( \frac{\lambda}{\lambda+\mu}\right) \ for \  i \geq 0,
	\end{equation*} and 
	\begin{equation*}
		P(i\ customers \ in \ the \ orbit \ and \ server\  is \ idle) = \left( 1-\frac{\lambda}{\theta}\right) \left( \frac{\lambda}{\theta}\right)^{i} \left( \frac{\mu}{\lambda+\mu}\right) \ for \ i \geq 0,
	\end{equation*}
	where $\lambda$ is the customer's arrival rate following the Poisson process, $\mu$ is the exponentially distributed service time, and the lead time for replenishment follows an exponential distribution with parameter $\gamma$.
	 The long-run performance measures were obtained using the joint distribution and the optimal pair $(s, S)$.

	 The article authored by Shajin \textit{et al.} \cite{skd} considered a QIS with reservation of inventory items, the possibility of cancellation, and over-booking for all future times up to a maximum of $K$ time frames. They analyzed the evolution of a system designed to reserve some items in advance by customers arriving at random moments. The arrivals of customers are according to a Markovian arrival process, and the service and inter-cancellation times are independent and exponentially distributed. They considered Poisson arrivals and exponentially distributed with rate $\lambda$, CLT (common lifetime) with rate $\gamma$, and  exponentially distributed amount of time with parameter $\mu$ as a special case and obtained a product-form solution that is for $\lambda < \mu$,
	\begin{equation*}
		x_{n}^{'} = \left(1-\frac{\lambda}{\mu}\right) \left(\frac{\lambda}{\mu}\right)^{n} \mbox{\boldmath{{$\xi^{'}$}}} \quad n \geq 0,
	\end{equation*}
	where  $\xi^{'}$ is the steady-state vector of the infinitesimal generator $\Tilde{B}$ corresponding to the system with negligible service time. Blocking of arrivals combined with the special case of Poisson arrival occurs when all time frames $1, 2, \dots, K$ are overbooked. They constructed a cost function and investigated its properties numerically.

	 Anilkumar and Jose \cite{mpank} discussed a PIS with positive service time and $(s, S)$ order policy. Wherein they assumed that the accessible inventory level is less than a prefixed level, it could minimize the loss of customers due to the stock-out period by decreasing or increasing the production rate. They investigated the inventory cycle using a stochastic decomposed solution for the steady-state probability vector. They obtained an explicit solution for the relevant performance measures and the closed-form solution to the model. Customer arrival follows a Bernoulli process with parameter $p$, and service time follows a geometric distribution with parameter $q$. The steady-state probability vector of the discrete-time Markov chain is, $\prod = (\pi_{0}, \pi_{1}, \pi_{2}, \dots) $ is,
	\begin{equation*}
		\pi_{i} =\begin{cases}
			\frac{q-p}{q} \mbox{\boldmath{{$\hat{\pi}$}}} &   i=0, \\
			\frac{q-p}{q} \frac{p}{pq} \rho^{i-1}  \mbox{\boldmath{{$\hat{\pi}$}}} 
			&  i\geq 1;
		\end{cases}
	\end{equation*}
	where $\rho = \frac{p\bar{q}}{\bar{p}q}$
	and $ \mbox{\boldmath{{$\hat{\pi}$}}} = (\hat{\pi}_{0},\hat{\pi}_{1}, \dots, \hat{\pi}_{s}, \hat{\pi}_{s+1,0}, \hat{\pi}_{s+1,1}, \dots, \hat{\pi}_{S-1,0}, \hat{\pi}_{S-1,1},\hat{\pi}_{S})$ is the steady-state probability vector.

	Alnowibet \textit{et al.} \cite{akaa} introduced a new stochastic mathematical model for inventory systems with lead times and impatient customers under deterministic and uniform order sizes. The model identified the performance measures in a stochastic environment, analyzed the properties of the inventory system with stochastic/probabilistic parameters, and validated the model’s accuracy. Balance equations were derived from a mathematical characterization dependent on the Markov chain formalism to analyze the system and obtained the steady-state probability as,
	\begin{equation*}
		P_n = \begin{cases}
			\dfrac{\left(\dfrac{\lambda}{\mu}\right)^{n}}{\left[{1+\dfrac{\lambda}{\mu}+
					\sum_{i=0}^{n-1}
					\dfrac{\lambda^{i+1}}{\prod_{j=1}^{n-1}\mu(\mu+j
						\alpha)}}\right] } & 0\leq n < N, \\
			0, & otherwise,
		\end{cases}
	\end{equation*}
	where $\lambda$ is the Poisson arrival rate, $\mu$ is an exponentially distributed service time, and the customer abandons the queue with rate $\alpha$. They achieved the performance by examining the graphical representation of the service process in steady-state as a function of both arrival distribution and the customer patience coefficient.

    Daduna \cite{dad} proposed a network of QIS where a common central server replenishes the inventories. The author designed a model for integrating the characteristics of manufacturing at various locations, inventory holding and control, replenishment production, and transportation with dispatching by system state-dependent adaptive regimes. Consumers behave in a way that combines back-ordering and lost-sales principles when they arrive at the locations. The production represents a collection of $J$ exponential single server queues connected to a closed network representing the production's support. He developed two models that allowed him to explicitly compute the stationary distribution of each system in product-form as well as solve the integrated location problem. That is Daduna construct a Markov process $Z = \{Z(t): t\geq 0\} $, If $Z$ is ergodic its stationary distribution $\pi$ is with normalization constant $G(b_{1}, b_{2}, \dots, b_{j}; \Bar{J})$  as,

    \begin{center}
     \begin{eqnarray*}
     \pi (m_{1}, k_{1}, n_{1}; \dots m_{j}, k_{j}, n_{j}; \dots m_{J}, k_{J}, n_{J})  = \prod_{l=1}^{J} \left\{ \left( \prod_{i=0}^{m_{l}+ k_{l}-1}  h_{l}(i)\right)   \left( \prod_{i=0}^{m_{l}} \frac{\nu d_{l}(x)}{i} \right)  \left( \frac{\nu}{\lambda_{l}}\right)^{k_{l}} \left( \prod_{i=0}^{n_{l}} \frac{\lambda_{l}}{\mu_{l}(i)}\right) \right\}  . \\ \left(\prod_{i=0}^{m_{1}+k_{1}+\dots m_{J}+k_{J}-1} h(i) \right) .  G(b_{1}, b_{2}, \dots, b_{j}; \Bar{J})^{-1}.
    \end{eqnarray*}   
    \end{center}
    If $Z$ is ergodic, the normalization constant is,
    \begin{eqnarray*}
        G(b_{1}, b_{2}, \dots, b_{j}; \Bar{J}) = \sum_{n_{1}=0}^{\infty} \dots \sum_{n_{J}=0}^{\infty} \left(\prod_{l=1}^{J} \left( \prod_{i=1}^{l} \frac{\lambda_{l}}{\mu_{l}(i)}\right) \right) . H(b_{1}, b_{2}, \dots, b_{j}; \Bar{J}),
    \end{eqnarray*}
    with 
    \begin{center}
      \begin{eqnarray*}
         H(b_{1}, b_{2}, \dots, b_{j}; \Bar{J}) = \sum_{(m_{1}, k_{1}; \dots m_{J}, k_{J}): 0 \leq m_{j}+k_{j}\leq b_{j}, j\in \Bar{J}}  \prod_{l=1}^{J} 
         \bigg\{ \left( \prod_{i=0}^{m_{l}+ k_{l}-1}  h_{l}(i)\right)  \\ \left( \prod_{i=0}^{m_{l}} \frac{\nu d_{l}(x)}{i} \right) \left( \frac{\nu}{\lambda_{l}}\right)^{k_{l}} \left( \prod_{i=0}^{n_{l}} \frac{\lambda_{l}}{\mu_{l}(i)}\right) \bigg\} . \left(\prod_{i=0}^{m_{1}+k_{1}+\dots m_{J}+k_{J}-1} h(i) \right),
    \end{eqnarray*}  
    \end{center}
    where $\Bar{J} = (1,2, \dots, J)$ is the set of locations, $\lambda_{j}$ is the Poisson arrival rate at location $j$, $\mu_{j}$ is the exponentially distributed service rate at location $j$ and $d_{l}(x) $ is the distance between $J+1$ and $l$.
    The most significant advantage of this model is the ability to express the decision for the location of the central manufacturing unit as a generalized Weber issue.

    Manikandan \cite{mr} considered a single server PIS with $(s, S)$ order policy. The author discussed a few realistic assumptions in the proposed model. For example, the produced item is not always necessary to add to the inventory due to manufacturing defects; and the items need not be served to the customer after the service completion epoch due to several reasons. More precisely, the item to be added to the inventory with a probability $\delta$. An item from the inventory is served to a customer with probability $\gamma$  at the end of his/her service, and the customer leaves the system empty-handed with $1 - \gamma$. The production takes a random time and follows an exponential distribution with parameter $\beta$. A manufacturing defect does not always add a produced item to the inventory: it is accepted and rejected with probability $\delta$ and $1-\delta$, respectively. Under these assumptions, the author derived the system's steady-state equation as a stochastic decomposition of the marginal probability distribution of a number of customers/demands and the number of items in the inventory. That is, under the necessary and sufficient condition $\lambda < \mu$, 
    \begin{eqnarray*}
    x_{i} = (1 - \rho) \rho^{i} \mbox{\boldmath{{$\pi$}}}, \quad  i \geq 0,
    \end{eqnarray*}
    where $\rho = \frac{\lambda}{\mu}$ and \mbox{\boldmath{{$\pi$}}} is the inventory-level probability vector, 
    \begin{eqnarray*}
        \pi(s-j) = \pi (S) \frac{\gamma \lambda}{\delta \beta - \gamma \lambda} \left( 1- \left(\frac{\gamma \lambda}{\delta \beta}\right)^{S-s}\right) \left(\frac{\gamma \lambda}{\delta \beta}\right)^{j}, \quad 0 \leq j \leq s;
    \end{eqnarray*}

    \begin{eqnarray*}
        \pi (j, 0) = \pi (S), \quad s + 1 \leq j \leq S - 1;
    \end{eqnarray*}

     \begin{eqnarray*}
        \pi(s-j) = \pi (S) \frac{\gamma \lambda}{\delta \beta - \gamma \lambda} \left( 1- \left(\frac{\gamma \lambda}{\delta \beta}\right)^{S-j}\right),  \quad s+1 \leq j \leq S-1,
    \end{eqnarray*}
    where
    \begin{eqnarray*}
        \pi (S) = \frac{\left(\frac{\gamma \lambda}{\delta \beta}-1\right)^{2}}{\left(\frac{\gamma \lambda}{\delta \beta}\right)^{S+2} - \left(\frac{\gamma \lambda}{\delta \beta}\right)^{s+2} - (S-s)\left(\frac{\gamma \lambda}{\delta \beta}-1\right)}.
    \end{eqnarray*}
    Derived the expected length of a production cycle and presented a few level-crossing results. Also, he provided a comparison of the performance measures for a few $(\gamma,\delta)$ pairs.

	Another recent article by Otten  \cite{os} studied a supply chain consisting of PIS at several locations, which a common supplier couples. Where item routing depends on the inventory at hand to obtain "load balancing." Under a constant review base stock policy, the supplier produces raw materials to replenish local inventories. The service starts immediately if the server is ready to serve a customer ahead of the line and the inventory has not been exhausted. Otherwise, the service begins when the next replenishment arrives at the local inventory. The author achieved to obtain the stationary probability distribution as a product of the marginal distributions of the production and inventory-replenishment subsystem.
	The unique limiting and stationary distribution is of the form, 
	\begin{equation*}
		\pi (\textbf{n},\textbf{k}) = \xi(\textbf{n}) \cdot \theta (\textbf{k}),
	\end{equation*}
	with \begin{equation*}
		\xi(\textbf{n}) = \prod_{j\in \bar{J}} \xi_{j}(n_{j}),\quad \theta(k) = \lim_{t \to \infty} P((Y_{1}(t),\dots,Y_{j}(t),W_{j+1}(t))= k);
	\end{equation*}
	\begin{equation*} \xi_{j}(n_{j}) = C_{j}^{-1} \prod _{l=1}^{n_{j}} \frac{\lambda_{j}}{\mu_{j}(l)}, \quad n_{j}\in N_{0}, j\in \bar{J}, 
	\end{equation*} 
	the normalization constants $C_{j}$ is,
	\begin{equation*}
		C_{j} = \sum_{n_{j}\in N_{0}} \prod_{l=1}^{n_j} \frac{\lambda_{j}}{\mu_{j}(l)},
	\end{equation*}
	where $\lambda_{j}$ is the customer arrival rate according to Poisson process and $\mu_{j(l)}$ is the service intensity.
	They derived an explicit solution for the marginal distribution of the production subsystem. Also, they derived an explicit solution for some special cases for the marginal distribution of the inventory-replenishment subsystem.

     A QIS with a random order size policy and lost sales, where the server takes multiple vacation policies when the on-hand inventory depletes, is considered by Zhang  \textit{et al.} \cite{zy}. The authors assumed that the customers prevent from entering the system when the on-hand inventory level is zero or when the server is off due to a vacation. Order size decisions and lead times are independent of the arrival process, the customer service time, and the server's vacation. The customer's arrivals are Poisson with rate $\lambda$, exponentially distributed service times of rate $\mu$, and $q_{i}$ is the probability that the replenishment order size is at least $i$. The authors derived the stationary joint probability distribution of the queue length, the on-hand inventory level, and the server's status in explicit product-form. And they also obtained that the stability condition was independent of the vacation rate $\theta$, the lead time parameter $\eta$, and the random order size distribution $\bar{p}$.
	 Thus, the steady-state probability vector of the process for if $\rho = \frac{\lambda}{\mu}<1$ as,
	 \begin{equation*}
		\mbox{\boldmath{{$x$}}}=(\mbox{\boldmath{{$x_{0}$}}},\mbox{\boldmath{{$x_{1}$}}},\dots),\quad\text{where}, \ \mbox{\boldmath{{$x_{n}$}}} = (1-\rho)\rho^{n} \mbox{\boldmath{{$\hat{\pi}$}}}, \quad n\geq 0,
	 \end{equation*}
	 where $\mbox{\boldmath{{$\hat{\pi}$}}} = (\hat{\pi}(0,0), \hat{\pi}(1,0), \hat{\pi}(1,1), \dots, \hat{\pi}(M,0), \hat{\pi}(M,1))$ is the steady-state probability vector of the generator $\{\hat{S}(t), t\geq 0\} $ and 
	 \begin{equation*}
	     \hat{\pi}(0,0) = \frac{\lambda}{\eta} K_{v}^{-1}, 
	     \end{equation*}
	 \begin{equation*}
	     \hat{\pi}(i,0) = \frac{\lambda}{\theta} p_{i} K_{v}^{-1}, \ i=1,2, \dots, M,
	 \end{equation*}
	 
	 \begin{equation*}
	    \hat{\pi}(i,1) = q_{i} K_{v}^{-1}, \ i= 1,2, \dots, M, 
	 \end{equation*} where
	 \begin{equation*}
	    K_{v} = \frac{\lambda}{\theta} + \frac{\lambda}{\eta} + \bar{p}. 
	 \end{equation*}
	 From the above relations, we can see that the stationary distribution of the system has a product-form of two marginal distributions. One is the stationary queue length distribution in the $M/M/1$ traditional queueing system, and the other one is the stationary distribution of the on-hand inventory level of the QIS system with the server's multiple vacations and random order size policy when the service time is negligible.
	 They also derived conditional distributions of the on-hand inventory level when the server is off due to a vacation or depleted inventory and when the server is on and working. The server's vacation provides several necessary system performance measures, and they analytically investigate the effect of the server's vacation on the performance measures and present some numerical results. This work gives a model simulation study with general arrival processes and service time distributions.

     Linhong Li \textit{et al.} \cite{llwz} studied a QIS with batch demands and multiple vacations. In this model, they considered a random order size policy for replenishment. The arrival of customers follow the Poisson process with rate $\lambda$ and exponentially distributed service time with rate $\mu$. When the on-hand inventory level reaches zero, the system sends a replenishment order and takes multiple vacations simultaneously. Customers arriving during this period are assumed to be lost. Customers joining the system during a vacation may have to wait a long time due to the randomness of the replenishment order size and the demand size, which reduces service satisfaction. The server will go on another vacation once it has finished it's current one and discovered the system is empty to solve the problem. They obtained stationary joint probability in product-form, which is given by,
    \begin{eqnarray*}
    x(n) = \left(1 - \dfrac{\lambda}{\mu}\right)\left(\dfrac{\lambda}{\mu}\right)^{n} \mbox{\boldmath{{$\hat{x}$}}}, \quad n \in \mathbb{N},
    \end{eqnarray*}
    where \mbox{\boldmath{{$\hat{x}$}}} is the stationary probability vector, $\mbox{\boldmath{{$\hat{x}$}}} =  (\hat{x}(0, 0), \hat{x}(1, 0), \hat{x}(1, 1), \dots , \hat{x}(S, 0), \hat{x}(S, 1))$ and they obtained a numerical solution for \mbox{\boldmath{{$\hat{x}$}}} using an algorithm. Important system performance measures were derived and numerically illustrated the impact of performance measures.

\subsection{\textbf{Single Server QIS}}
In this section, we focus on the contributions of a single-server QIS. A single-server QIS is a system in which customers arrive at a single-server queueing system in which the server uses some commodities from the inventory to fulfil the service.

 A. Melikov \textit{et al.} \cite{meli} investigated single-server queuing-inventory system (QIS) models with catastrophes in the warehouse and negative customers (n-customers) in the service facility. Customers (c-customers) who came to purchase inventory could be queued in an infinite buffer. Customers in the system (on a server or in a buffer) are still waiting for stock replenishment even after the system's entire inventory is destroyed via a catastrophe. They considered two replenishment policies (s, S) policy or randomized replenishment policy. They obtain the stability conditions for the systems under investigation and demonstrate that, in particular situations, they are independent of the size of storage, the rate of catastrophes, and the rate of replenishment. They developed formulas that allow analyzing of the effect of the initial parameters on performance measures of the studied QISs as well as on expected total cost (ETC) and appropriately select the optimal RPs parameters so that the ETC is minimized.

	Chakravarthy \textit{et al.} \cite{csr} analyzed a single-server queueing-inventory model in which the customers served in varying-sized batches based on the predefined thresholds and the items available in the inventory. They used an  $(s, S)$ control policy for the proposed model. The model was analyzed using the Matrix-Analytic Method $(MAM)$ at first, with all underlying random variables assumed to be exponentially distributed, and further discussed an outline of the model in a more general setup. Because of the complexity of the model, the authors adopted the simulation method for good validation in the case of the analytic counterpart of the exponential model when they made more general assumptions on the underlying random variables. Several numerical examples were provided to accomplish their analysis.

	An $(s, S)$ inventory policy with exponential replenishing time is considered in Amirthakodi and Sivakumar \cite{ams}. It provides a busy period analysis of the system and the waiting time distribution for both primary and feedback customers using the Laplace-Stieltjes transform. The authors used an algorithmic approach to provide numerical illustrations for the convexity of the long-run total expected cost. This paper significantly extends to include feedback from customers to multi-server inventory systems.

	Anilkumar and Jose \cite{amp} investigated a discrete-time $(s, S)$ single server queueing inventory model with self-interrupted priority, Bernoulli arrivals, and geometric distributions for service time and lead time. The authors proposed two types of queues: high priority and low priority, with the arriving customers entering into a high-priority queue of infinite capacity. In addition, if the service is interrupted, then the interrupted customer will be forced to move to a lower-priority queue of infinite capacity. An item in the inventory is provided to the customer and leaves the high-priority queue at the service interruption/completion epoch. Customers in the lower priority queue get service by preemptive priority discipline. The system is analyzed using $MAM$. They obtained marginal distributions of queue length for both primary and interrupted customers.

	Maqbali \textit{et al.} \cite{Maqb}  discussed an $M/PH/1$ queueing system with an inventory level under $(s, S)$ policy motivated by the hospital inventory management model. The arrival of customers follows Poisson Process with the rate $\lambda$, and the service time attached with inventory level follows $PH$-distribution with representation $(\beta, T)$ of order $m$. The order decision is made based on a reorder point $s$, and the lead time determines the replenishment, which follows an exponential distribution with the rate $theta$. The authors concluded that the expected number of customers in the system $E[N]$ and $E[N_0]$, due to lack of inventory decrease as $S$ increases numerically. Besides this, the expected number of items in inventory $E[I]$ increases when $S$ increases. Due to a shortage of inventory $b_1$, the probability that a customer waits for service decreases, and the probability that the server is idle $b_0$ has almost no change as $S$ increases.

	R Manikandan and S S Nair \cite{rms} studied an $M/M/1$ QIS under $(s, Q)$ replenishment strategy with server working vacations, vacation interruptions, and lost sales. The authors assumed that the server takes vacation only in the absence of customers in the system and not due to the inventory level falling to zero at a service completion epoch. The server provides service at a lower rate during working vacations than in the normal service mode. With the system having infinite capacity, they obtained the stability condition for the system. They performed the busy period analysis and derived the stationary waiting time distribution in the queue and also evaluated various performance measures. Numerical illustrations are given to show the impact of the system's performance and discussed an optimization problem.

	Krishnamoorthy \textit{et al.} \cite{ak}, considered a single-server QIS with $BMAP$, and subsequent arrival batch sizes form a finite first-order Markov chain. It provides a control policy in which the server goes on vacation to assure idle-time utilization when a service process is frozen until a quorum may commence the next batch service. During the vacation, the server generates inventory for future services until it reaches a predetermined level or the number of customers in the system reaches a maximum service batch size, or whichever comes first. The time taken to process one unit of inventory follows a phase-type (PH) distribution. They estimated the steady-state probability vector of the infinite system in this paper. The distributions of inventory processing time, idle time, and vacation cycle length in a vacation cycle are determined. Using numerical examples, they compared the proposed system to a queuing-inventory system without the Markov-dependent assumption on successive arrivals and service batch sizes.

	Dudin and Klimenok  \cite{dakv} studied a wireless sensor network node with an energy harvesting model as a single server $MAP/G/1$ QIS. The authors assumed that an arriving customer will get service only if an energy unit is available. Also, if the server is idle or it has no energy units available, the system places the arriving customers in the infinite-capacity buffer. They examined the stationary behavior of the sojourn time and the process of system states at any moment. The system under investigation only allows a customer to receive service if it draws a unit of energy from a stock of finite size. They obtained the steady-state probabilities, and a condition for the existence of the stationary distributions of these processes was derived. In addition to the Laplace-Stieltjes transform of the distributions of the virtual and actual sojourn times, formulas for performance measures are derived. Illustrating numerical examples highlights the significance of accounting for input flow correlation in understanding how the system performance measures behave with respect to its parameters.

	\subsection{\textbf{Multi-Server QIS}}\label{subsec3}
	
	A multi-server queue has two or more service facilities simultaneously providing identical service. When the servers use some commodities from the inventory to fulfil the service, the system becomes a multi-server QIS.

	A multi-server QIS with non-homogeneous Poisson arrivals and exponentially distributed service times and lead times are considered Yue D \textit{et al.} \cite{yzy} in 2016. The readers may note that this paper is nowhere discussed in any of the previous survey articles on QIS. Wherein the authors discussed a continuous review $(s, S)$ policy and formulated the system as a quasi-birth-and-death (QBD) process, obtaining stability conditions for the system. When the on-hand inventory level reaches the reorder point $s$, a variable replenishment quantity is placed to restore the on-hand inventory level to $S$ upon replenishment. They also considered an emergency replenishment policy, in which one item is replenished immediately from a second source with no lead time when the on-hand inventory level drops to zero. The stationary condition of the system, the stationary distribution of the joint queue length, and the on-hand inventory level are derived using a matrix geometric method. Finally, the authors formulated a cost function to minimize the optimal cost, and to obtain the optimal $(s, S)$ replenishment policy.

	Suganya \textit{et al}. \cite{cs} analyzed an $(s, S)$ inventory system with two heterogeneous servers and multiple vacations. They assumed arrivals to be $MAP$ and two parallel servers that offer customers heterogeneous phase-type services. The most significant contribution of their study is to allow vacations/multiple vacations for each server. When it comes to inventory, the vacation begins when the customer level drops to zero and when the inventory runs out. Moreover, even if there is at least one customer in the system at the end of the vacation and there are no items in stock, another vacation will begin. The model has broader applications since it includes these real-time characteristics. The authors obtained the stationary distribution of the system state, the joint probability distribution of the number of customers in the system, the inventory level, and the server status in the steady-state are derived.

	Benny \textit{et al.}  \cite{bbv} investigated a multi-server queueing-inventory model in which one type of customer is encouraged to serve another type of customer, which improves the efficiency of the service facility and is helpful in crowd-sourcing. They assumed that the items in the inventory to be served is finite. Thus, when the item is not available service cannot be provided. The authors numerically illustrated the effect of the revenue function and  probability on the revenue function for single, two, and three servers.

	Hanukov \textit{et al.} \cite{Hanv} proposed an $M/M/2$ model in which idle servers produce and store preliminary services to reduce the sojourn time of incoming customers and increase customer arrival rate. They formulated the process as a three-dimensional continuous-time Markov chain and obtained closed-form solutions using matrix geometric analysis. Analytically it is proved that the stability condition of this model is identical to that of a standard $M/M/2$ queue. An optimization problem to maximize the server's expected profit by controlling the capacity of preliminary services and the investment in increasing customers' arrival rates when preliminary services are available. The numerical examples provided in the article revealed that the investigated profit function is quasi-concave.

 Jegannathan \textit{et al.} \cite{jeg} investigated an intricate multiserver queueing-inventory system that incorporates an asynchronous server vacation and customer retrial facility. The system is characterized by c identical servers, a finite-size waiting space, and an S-item storage area. Service times are modeled with an exponential distribution. After the busy period, if each server encounters insufficient customers or items in the system, they start its vacation. When the server's vacation is over and it realizes there won't be a chance of getting busy, it enters an idle state if there aren't enough customers or items; otherwise, it will take another vacation. The vacation time for each server also occurs independently from that of the others. For inventory management, the system adopts an $(s, Q)$ control policy for replenishment.  The Neuts and Rao matrix geometric approximation approach, is used in this model to determine the stability condition and stationary probability vector.

	Rasmi \textit{et al.} \cite{rk}, discussed a QIS with heterogeneous customers of $K$ types arriving according to a Marked Markovian Arrival Process $(MMAP)$. The system classifies customers into different classes based on the service they seek and assigns different priorities to each class, resulting in different levels of inventory admitted to exhaust. Each class has a single service node that provides exponential services with class-dependent service rates. A single source of inventory replenished according to an $(s, S)$ policy with exponential lead times serves customers of all classes. Using $MAM$, they derived the steady-state conditions and invariant probability measure. A detailed analysis of inventory recycle time; and presented an optimization problem.

	Shajin \textit{et al.} \cite{skms} considered a multi-server PIS with an unlimited waiting line. The arrivals of customers follow a non-homogeneous Poisson process, and the production time for a unit item is phase-type distributed. They also investigated the model with exponentially distributed production time as a special case. An emergency replenishment policy is adopted to prevent the loss of demand when the inventory level is equal to zero. One item is replenished immediately with zero lead time when the on-hand inventory depletes to zero. The emergency replenishment policy ensures no decrease in demand in this model. When the on-hand inventory reaches zero, a local purchase (emergency replenishment) brings the inventory level to 1. The impact of various performance measures and an optimal cost function is demonstrated.

 \subsection{ \textbf{QIS with Heterogeneous Customers}} \label{subsec4}

  Modeling and analyzing queueing systems with heterogeneous customers is an important area of research in queueing theory. Heterogeneous customers refer to a scenario where multiple types or classes of customers are in a queueing system that differs from homogeneous customers. Each customer class may vary in terms of arrival rates, service requirements, and priorities.

	Baek \textit{et al.} \cite{jb}  analyzed a single-server QIS with a $MMAP$ for heterogeneous customers. They considered preemption priority for type-1 customers to be limited compared to type-2 customers. Type-2 customers have an infinite buffer, while type-1 customers have no buffer. Depending on the type of customer, a fixed number of additional consumable items are required for service. They proposed a threshold technique for admitting type-2 customers to service to minimize type-1 customer loss probability while maximizing system throughput. Only when the server is idle, and the quantity of additional consumable products in stock exceeds the fixed threshold type-2 customer then the service begin. They computed the waiting time and stationary distributions of the system states, provided a numerical example that highlights some interesting effects of performance measures, and demonstrated the optimization problem.

	Priority multi-server retrial QISs with Markovian arrival processes generated by a single or dual source and exponentially distributed service times, with finite queueing and orbit spaces, were investigated by Wang \textit{et al.} \cite{wfb}. To handle the non-identical services, the authors proposed a generalized stochastic Petri net (GSPN) model. They assumed that the Type 1 customer has higher priority than the Type 2 customer. An arriving customer of either type 1 or type 2 is taken for service immediately when there are sufficient resources; at least one server and one inventory are available simultaneously. If resources are unavailable, type 1 customers enter into a waiting line, or if the waiting line is full, loss occurs.
    Similarly, when resources are unavailable, type 2 customers enter into an orbit or get lost when the orbit buffer is full. The authors discovered that the performance is unaffected by different $MAP$ input generations. They found that performance can vary significantly if service rates are different. Therefore, it is useful to have an analytical model that can handle service heterogeneity, such as the GSPN model. To demonstrate the applicability of the proposed GSPN models, they solved an optimization problem to determine the optimal threshold for state-dependent servicing protocol.

	Shajin \textit{et al}. \cite{sdjvvm} considered a queueing-inventory problem arising in transporting passengers in which seats in the passenger's vessel assume physically available inventory. The arrival of customers forms a marked Poisson process and exponentially distributed service time. Two classes of customers (Low-priority and High-priority) arrive at the system to reserve seats. The system permits over-booking (if no item is available in the system, a limited number of reservations can be made; this is known as overbooking) by setting an upper bound. The vessel capacity for the scheduled departure is adjusted based on the number of over-bookings. Low-priority customers join an infinite capacity queue transfer to the buffer for service. High-priority customers are given non-preemptive priority over low-priority customers. Customers of both types have service times that follow independent exponential distributions with parameters based on the stage of the common lifetime. The \textit{common lifetime} is a term used to describe the potential expiry date of each item on hand; in some cases, it may be the same lifetime for all. Finally, they discussed an optimization problem to find the optimal cost.

	A continuous review inventory $(s, Q)$ finite buffer queuing model with two heterogeneous servers and mixed priority service was analyzed by Kingsly \textit{et al.} \cite{sjk}. In this model, one of the two servers is for high-priority customers, and the other serves both customers with a mixed-priority service scheme. High-priority customers actually receive service after selecting the purchased item from the inventory, while low-priority customers only arrive for repair work. If there is no high-priority or low-priority customer in the buffer, or if there is only one inventory in the system, the flexible server takes multiple vacations. The authors obtained the steady-state joint distributions of the waiting time of mixed priority service customers and the inventory level by matrix method.

	A two-server QIS consisting of two-parallel queues and two categories of customers is reported by Jeganathan \textit{et al.} \cite{kjna}. High-priority customers generate demand for an item, whereas low-priority customers arrive only for repair. Server 1 offers service to the first queue involving only high-priority customers, and server 2 is a flexible server that offers its service to both queues. The authors considered that server 1 would be idle and server 2 will be on vacation when there is no service available in both queues. To minimize the loss probability of low-priority customers, they found a priority service rule and served the low-priority customers on threshold-based service. The authors identified significant changes in system measures by computing the limiting probability distribution of the four random variables in a steady state.

	Fathi \textit{et al.}. \cite{fmkm} investigated a location-inventory model with stochastic demand and multiple customer priority classes. This study conducts the $(s, Q)$ inventory control policy at the distribution centres (DCs) level. The model aims to find the optimal locations for DCs and their inventory policies simultaneously. The stock level of DCs is modeled as a Markov chain process and analyzed in the first phase. In contrast, in the second phase, a mathematical program is used to determine the number and locations of DCs, the assignment of retailers to DCs, and the order quantity and safety stock level at DCs. NP-hard, a hybrid genetic algorithm, was developed to solve the problem.

	The queueing–inventory model studied by Jeganathan \textit{et al.}  \cite{jkvsh} with two classes of customers, two stages of service for low-priority customers and a single stage of service for high-priority customers. Adopting the mixed-priority discipline, assumed that the arriving high-priority customer could interrupt the service of a low-priority customer whenever the inventory level is less than the threshold level during the stage-I service of a low-priority customer. Otherwise, it used non-preemptive priority discipline in both stages. The interrupted low-priority customer moves to orbit and retries for the service whenever the server is free. The inventory replenishment follows the $(s, Q)$ ordering policy, with the exponentially distributed lifetimes of the items. An expression for the stability condition is determined explicitly. Numerical examples for different sets of input values of the parameters are computed to illustrate the system characteristics.

	\subsection{\textbf{Multi -Commodity QIS}}\label{subsec5}
	
	A QIS involving two or more distinct commodities is said to be a Multi-Commodity QIS. There are only a few works reported in this area. Those are listed as follows:

	In a queueing-inventory model, customers arrive according to a point process; each customer demands one or more items, but not exceeding a predetermined (finite) value is considered by Chakravarthy \cite{chasr}. The authors analyzed two models: In Model 1, any arriving customer finding the inventory level to be zero is considered lost, whereas in Model 2, the loss of customers occurs in two ways. First, an arriving customer finding the inventory level to be zero with the server being idle will be lost. Secondly, the customers present at a service completion with zero inventory will all be lost. Replenishment is based on the $(s, S)$ policy. In both models, the demands were partially met based on the requests and the availability of the items in the inventory. The steady-state analysis of the model is performed using $MAM$. This study evaluated and analyzed two models that deal with batch demands to demonstrate the practical benefits of each.

	Serife Ozkar \textit{et al.} \cite{so} examined a two-commodity $(s, S)$ QIS with two classes of customers in which customers request random amounts of service. The priority customers receive a non-preemptive priority over other customers. The model is constructed as a quasi-birth-and-death structure. The system performance measures and steady-state distributions are obtained. Also, provided the closed form of the system load. As in previous studies of single commodity-single customer inventory systems, the stability of the defined system does not depend on inventory parameters. Instead, the system load is determined by the arrival rate of priority customers. Additionally, they found that the system load is affected by the buffer size, and the system can use a smaller buffer size if the arrival rate is high.

	A multi-commodity QIS with one essential and $m$ optional items under the $(s, S)$ ordering policy with lead time is analyzed by Shajin \textit{et al.} \cite{sd}. The authors introduced the concept of optional items for service in this paper and assumed that the service is not provided without an essential item. The customer either leaves the system immediately after the service of an essential item or he/she goes for an optional item(s). The customer can demand more than one optional item, and if the demanded optional item is not available, the customer leaves the system after purchasing the essential item. A single server system with $MAP$ arrivals, $PH$-type  service time distribution, and exponentially distributed optional items are considered. Also derived several performance measures and did optimization problems.

	Jeganathan \textit{et al.} \cite{jkrma} considered a single server two commodity inventory system with queue-dependent services for finite queue and an optional retrial facility. Here are the demands for two commodities in a Markovian inventory system, where one item was designated as a major item and the other as a complimentary item. The authors investigated a strategy of $(s, Q)$ type control for the major item with a random lead time but instantaneous replenishment for the complimentary item and used the classical retrial policy for the orbiting customers. They derived the joint probability distributions for commodities, the number of demands in the queue, and the orbit in the steady state. The effect of the various performance measures is discussed  using numerical examples which are associated with diverse stochastic behaviours.

	The performance of a perishable QIS under $(s, Q)$ policy for two commodities with optional customer demands were analyzed by Anbazhagan \textit{et al.} \cite{anjgp}. The authors assumed that all customers who come to the system can only purchase the first or second item or service. The purpose of this work is to show the significance of the impact of optional demands on the system’s performance. Customer arrivals are assumed to be the Markovian Arrival Process ($MAP$). Using matrix-geometric methods, computed the joint probability distribution for the first commodity, the second commodity, and the number of customers in the system in the steady-state case. Various system performance measures are derived and provide a numerical illustration of the optimal value for diverse system parameters.

	Ozkar S \cite{oz} studied a two-commodity QIS with two types of customers. Type 1 customers have a non-preemptive priority and a finite waiting space over Type 2 customers with no limit in waiting room capacity. They considered two inventory policies separately for the two commodities. They added a local purchase option of order size one with instantaneous replenishment to avoid losing the customer waiting in the queue. The quasi-birth and death structure is constructed, and performed the steady-state analysis using the matrix-geometric method. The sensitivity analysis of the parameters, the optimum buffer size, the optimum service rate for different values of parameters, and the optimum inventory policy were obtained for each commodity.

	\subsection{\textbf{Queueing Systems with Perishable Inventory}}\label{subsec6}
	
	Perishable inventory consists of products that expire or lose value over time, eventually becoming worthless. Now, this section is devoted to QIS with perishable inventory.

	Melikov and Shahmaliyev studied a QIS model with perishable inventory \cite{measa}. In this model, when the inventory level is zero, it is assumed that some demands do not acquire the item after service completion and that some demands are impatient in the queue. The inventory replenishment policy belongs to the $(S-1, S)$. The exact and approximate formulas for the calculation of the steady-state probabilities of the system are obtained. The numerical experiments they have given show the formulas' high accuracy. Finally, they have given an optimization problem of choosing the optimal server for cost minimization.

	G Hanukov \textit{et al.}  \cite{gb} considered a typical fast-food service system as a customer queueing system combined with an inventory of perishable products. The authors simultaneously applied time management policies and inventory management techniques, which helps to improve the efficiency of the system. Based on a combined queueing and inventory model, they proposed a method in which each customer's service was divided into two independent stages. The first stage is to be generic and complete without customers, whereas the second stage requires the presence of customers. When the system is empty of customers, the server produces an inventory of first-stage services. The steady-state probabilities are derived and formulated using the matrix-geometric method. Also showed that the production rate of preliminary services does not affect the stability of the system. An economic analysis is provided to explore the optimal preliminary service capacity and optimal level of investment in preservation technology.

	Koroliuk \textit{et al.} \cite{kv} analyzed the model of a QISs with a single server and perishable stock, in which impatient customers can create queues of varying lengths. They introduced a two-bin $(s, S)$ replenishment policy, in which a threshold $s, s< S$ was set, and if the system's inventory level is higher than this value, the system does not make replenishment orders and assumed that the stock that has already been distributed would not perish. Also, they suppose that the arrived customer will be lost with sure probability if there are $N$, $1< N < \infty$, customers in the system at the time instant,  also at the same time, any external customer arrival will be allowed into the queue. If there is no customer in the queue, the server takes multiple vacations with exponential vacation time. Customers are only impatient while in line; a customer on the server does not leave the system without being served. The server's status usually determines the degree of impatience of customers in the queue.  Accurate and approximate methods for determining the characteristics of the models are developed and used to conduct numerical experiments.

	Bhuvaneshwari  \cite{bm} discussed a continuous review perishable inventory system with a finite number of homogeneous sources of demands. The ordering policy is the $(s, S)$ policy. The lead time of ordering and the lifetime of each item are distributed exponentially. The authors considered that all the arriving customers demand the essential service first, and some of them may further demand one of the other optional services. The service times of the essential service and the optional services are assumed to be exponentially distributed. The joint probability distribution of the number of customers in the waiting hall and the inventory level for the steady-state case is obtained. And derived the Laplace-Stieljes transforms of waiting-time distribution of customers in the waiting hall, some important system performance measures in the steady-state, and the long-run total expected cost rate.

	Saranya and Shophia  \cite{ns} investigated a continuous review $(s, S)$ perishable inventory system with a finite number of homogeneous demand sources. The arrival of demands is a quasi-random process, and the perished items were replaced free of cost by the supplier at the time of replenishment. When the on-hand inventory level was zero, they supposed that any newly arriving demand was offered a choice of postponement, and the selection of postponed demands is based on some preset rules. The various system performance measures, the joint probability distribution of the on-hand inventory level, the number of perished items stored for replacement, and the number of customers in the pool in the steady-state case are obtained.

	Laxmi and Soujanya  \cite{lp} examined a continuous review perishable inventory system with two types of customers, positive and negative, arriving according to a $MAP$. The lifetime of an item and the reordering lead time is exponentially distributed. They assumed that demands occurring during a stock out or busy period are either lost or enter an orbit of finite size. Depending on the number of demands in the orbit, orbital demands compete for service at an exponential rate. After an exponentially distributed amount of time, the waiting demands in the orbit may renege the system. At zero inventory, the server takes multiple working vacations. The authors developed a recursive method to find the steady-state joint probability distribution of the number of customers in orbit and the inventory level. Finally, various performance measures and cost analyses are presented with numerical results. The model is appropriate for situations where the negative arrival removes a random number of ordinary customers from the system. The proposed model can be applied in various real-world scenarios, including telecommunication networks, ATM networks, etc.

	Melikov and Shahmaliyev  \cite{azmm} considered a perishable $(s, S)$  QIS with positive service time and different types of customers. The customer arrival follows the Markov-Modulated Poisson process with infinite queue length and finite inventory size. The items in the inventory are considered to perish, but the item reserved while serving the customer cannot perish. The customers in the queue become impatient when the inventory level drops to zero. If no items are left in the inventory, the arriving customer enters the queue or leaves the system according to the Bernoulli scheme. The stationary distribution of the system and performance measures are formulated using Gillespie's Direct simulation method. Using numerical experiments, they illustrated the dependence of performance measures on reorder level and the solution to the optimization problem.

	Jeganathan \textit{et al.}  \cite{jkma} focused on a perishable inventory model with two stations, two dedicated servers, and one flexible server, where service rates are different for stations and servers. The item could not perish while servicing. The joint stationary distribution of the number of customers in the system, the status of the three servers, and the inventory level were derived in steady-state using the $MAM$. They computed various system performance measures and calculated the long-run total expected cost rate. Extensive numerical illustrations were presented to show the effect of the change of values for the parameters that affect the total expected cost rate by assuming a suitable cost structure on the inventory system.

	A finite and infinite three-dimensional perishable QIS model with delayed feedback is proposed by Melikov \textit{et al.} \cite{ma}. The authors further assumed that the customers either leave the system with/without purchasing an item or join the orbit for the decision-making. The replenishment policy belongs to the $(s, S)$ class and is exponentially distributed with a lead time. They obtained the exact and approximate formulas for calculating the system's steady-state probabilities and performance measures, where the exact method is based on solving balance equations and is only suitable for finite perishable QIS. The approximate approach is based on the State Merging Algorithm (SMA) of Markov Chains and is applicable to finite and infinite systems. The limitations and advantages of the SMA algorithm are also explained. The authors presented the dependence of performance measures on the reorder levels of the finite model, and also they compared the same with the partly-infinite models and described the features of the model using graphical illustrations.

	Prasanna Lakshmi \textit{et al.} \cite{kpl} investigated a Markovian queuing-inventory model consisting of a single server with two-stage working vacations with the maximum stock capacity. The authors considered the service rate in the $1^{st}$ stage of working vacation to be lesser than the service rate in a regular busy period but greater than the service rate in the $2^{nd}$ stage of a working vacation. The replenishment policy of any order is $(s, Q)$ policy, and supposed that the inventory is perishable. The limiting distribution of all random activities in the system as stationary is derived using the $MAM$. The optimal behaviour of maximum inventory and reordered levels with various cost structures is discussed, and also studied the nature of the number of customers lost and the expected waiting time of a customer.
	
	Uzunoglu Kocer and Yalcin \cite{uk} examined a continuous review $(s, Q)$ inventory model for perishable items with two demand classes called priority and ordinary customers. When the on-hand inventory drops to a pre-specified level, they assumed that only the priority customer demands were met, whereas the demands from ordinary customers are lost. The demand occurring during stock-out periods is considered to be lost (it is also called lost sales). The steady-state probabilities for this model are obtained as, 
	\begin{equation*}
		P_j = \frac{\mu}{(\lambda_1+j\gamma)} A_j P_0, \quad 1\leq j\leq s,
	\end{equation*}
	\begin{equation*}
		P_j = \frac{\mu+\lambda_1+s\gamma}{(\lambda_1+\lambda_2+j\gamma)} \frac{\mu}{(\lambda_1+s\gamma)} A_s P_0, \quad s+1\leq j\leq Q,
	\end{equation*}
	\begin{equation*}
		P_j = \frac{\mu+\lambda_1+s\gamma}{(\lambda_1+\lambda_2+j\gamma)} \frac{\mu}{(\lambda_1+s\gamma)} A_s P_0 - \frac{\mu}{(\lambda_1+\lambda_2+j\gamma)} A_{j-Q}P_0, \quad Q+1\leq j\leq Q+s,
	\end{equation*}
	where  
	\begin{equation*}
		A_{j} = \begin{cases}
			1 & j\leq 1, \\ \prod_{n=1}^{j-1} ( \frac{\mu}{\lambda_1+s\gamma}+1) & j\geq 2,
		\end{cases}
	\end{equation*} and $ \lambda_1 $, $ \lambda_2$ are the demand rates of prioritized and ordinary customers according to the Poisson process, respectively, $ \mu$ is the exponentially distributed lead time and the lifetime of each item is exponentially distributed with rate $\gamma$. They derived the expected total cost function and proved the pseudo-convexity of the cost function. The numerical computations provided show that the perishability of the product critically affects the inventory levels and inventory costs accordingly, i.e., it shows the significance of perishable inventory management. The model can be applied in many areas as there are different types of customers who demand perishable products.

	A continuous-review perishable inventory model with random lead times and state-dependent Poisson demand was analyzed by Y Barron  \cite{bar}. The paper dealt with demand uncertainty and allowed for random batch demands. A comprehensive analysis of two main models with different lead times and perish times under back-orders or lost sales so that the models can be applied to many industries in situations where the system is subject to random perishability, random lead time, and demand uncertainty. The long-run average cost function under the $(S, s)$ replenishment policy was derived with a probabilistic approach. Numerical examples are used to demonstrate the impact of changing batch size and other system parameters on the optimal policy, which indicates that it performs better for a general perish time, although the Markovian policy can be used as a good approximation of the average total cost. It is shown that the optimal cost may differ for different average batch sizes, while the batch variability seems to provide some robustness.

	Barron Y and Baron O  \cite{by} considered a $(Q, r)$ perishable inventory system with state-dependent compound Poisson demands, random batch size, general lead times, exponential shelf times, and lost sales. The authors used Queueing and Markov Chain Decomposition (QMCD) approach to derive exact closed-form cost expressions. The authors provided  a substantial extension of QMCD and inventory theory for the Markov chain embedded at the arrivals of new batches. Also, a closed-form expression is obtained for the expected total long-run average cost function from which they obtained the optimal $Q$ and $r$ parameters. Using the demand size and lead time distributions the authors discussed the optimal control problem to minimize inventory management costs.

	Sangeetha and Sivakumar  \cite{sns} proposed the model of a finite capacity continuous-review perishable inventory system with a service facility of exponentially distributed service times and Markovian arrival. The $(s, S)$ replenishment policy is used, and the lead time follows the $PH$-type distribution. The authors assumed the items are perishable after a random time if kept on the shelf for a long time, with the lifetime of each item assumed to be exponential. The model is formulated as a semi-Markov decision problem. Using a linear programming algorithm, they computed the stationary optimal policy. They found that the total expected cost rate can be minimized by controlling the service rate at each instant of time; also, they inferred that the optimal policy depends on the arrival phases, replenishment phases, inventory levels, and the number of customers.

	A continuous review of perishable goods with lost sales and $(s, S)$ policy was considered by  Baron \textit{et al.}  \cite{bp}. The exponentially distributed perishability and the lead times are assumed, while two cases of demand distribution are considered, one is Poisson, and the second one is compound Poisson with general demand sizes. They derived the average cost per unit of time and the total cost comprised of four components known as ordering cost (both variable and fixed), holding costs, cost of unsatisfied demand, and cost of perishability. Under the average cost criterion, they obtained exact, closed-form expressions for the cost of the system under the average.

	A continuous review perishable inventory system working under the $(s, Q)$ policy with an infinite orbit and retrial policy analyzed by Reshmi and Jose \cite{rpsj}. The arrival of primary customers follows $MAP$, and the service time follows the $PH$-type distribution. Items are subject to decay at some linear rate. An unsatisfied primary customer enters an orbit or exits the system. The retrial customers try to access the idle server at a linear rate. The stationary probability vectors are obtained using the $MAM$ with different parameter values. The optimal profit function is computed numerically.

	Using the QMCD approach, Y Barron and O Baron  \cite{ybo} provided a detailed analysis of an $(S, s)$ continuous-review inventory model for perishable items with random lead times and times to perishability, and also when demand is state-dependent. The authors presented a comprehensive analysis of two main models, for the first model, they assumed a general random lifetime and an exponentially distributed lead time. In contrast, an exponentially distributed lifetime and a general lead time were assumed for the second model. Both models are analyzed under back ordering and lost sales assumptions. They considered a fixed cost for each order, a purchase cost, a holding cost, a cost for perished items, and a penalty cost in the case of shortage. The variability of lead time is more expensive than the variability in perishability time is established numerically.

	Radhamani et al. \cite{rvsba} considered a continuous perishable inventory system at a service facility with an infinite waiting capacity. The server takes a vacation during the stock-out period. During the stock-out period, the server can avail of vacation. Due to the lower reorder level than the inventory level, stocked items may perish during server vacation. Three replenishment policies are considered for the system: variable-size order policy, $(s, S)$ policy, and order up to $S$ policy, in which variable-size order policy is the best policy among the three in terms of total expected cost rate. A numerical study analyzes these policies by assuming $MAP$ for the demand process.

	\subsection{\textbf{QIS with Postponed Demands}} \label{subsec7}
	
	The demands during the stock-out period are either lost (lost sales) or satisfied only after the arrival of ordered items (backorders). In the case of back orders, the back-ordered demand may have to wait even after replenishing. This type of inventory problem is called inventory with postponed demands.

 A continuous review $(s, S)$ inventory system with a service facility consisting of a finite waiting hall and a single server is studied by Jenifer and Sivakumar  \cite{jsa}. The customers arrive according to a Poisson process and are served according to an FCFS queue discipline. An arriving customer who finds the waiting hall is full enters into the pool of infinite size or leaves the system, according to a Bernoulli trial. When no customer is in the pool, the system behaves as an $M/M/1/N$ queue maintained at a service facility. The pooled customers are selected from the pool on a specific basis. The system behaviour is described by the Quasi Birth-Death process. In the steady-state case, they obtained the joint probability distribution of the number of customers in the pool, the number of customers in the waiting hall, and the inventory level. Various stationary system performance measures and the total expected cost rate are calculated. The Markovian arrival process with a multi-server system is a challenging extension of the proposed model.

	The significance of modified and delayed working vacations on perishable QISs with two heterogeneous servers was discussed by Jeganathan and Abdul Reiyas \cite{kjm}. They supposed that one server is used for high-priority customers with modified working vacations who demand both items and services. In contrast, the other is for low-priority customers with delayed working vacations who demand service only. Items are replenished under the $(s, Q)$ ordering policy with a positive lead time. The server is partially involved in serving any customer on working vacation, and the service was completed with a positive rate lower than the normal service rate. In the meantime, the server can perform some other duties. The steady-state distribution of the system at the limiting case is obtained, and the distributions for waiting times of both customers are analysed. The authors inferred that any server's time on working vacation is much greater than that of the server in the normal service period in the long run. The authors demonstrate the advantages of the proposed model with the features of modified working vacations compared to the cases of simple vacations and non-delayed working vacations.

	Padmavathi \& Sivakumar \cite{ip} considered an $(s, S)$ inventory system with postponed demands and single server vacation under discrete time set up. The authors proposed three vacation policies and found the optimal and expected cost for all three policies. If the stock-out period occurs, the server becomes idle for a random period called the modified multiple vacation policy. The second vacation occurs whenever the inventory level is zero. The third vacation policy is considered a continuation of the first policy. The authors assumed discrete $MAP$ for arrivals and independent Discrete $PH$ for vacation, idle, and discrete time. They found that the multiple vacation policy is the best among the three policies. The authors numerically analyzed the impact of system costs on the total optimal cost of a modified multiple vacation policy.

	A continuous review $(s, S)$ inventory system with two types of customers arriving according to two independent Poisson processes is analyzed by Sangeetha \textit{et al.}  \cite{nst}. The authors' main objective is to find an optimal selection rate of pooled customers at any given inventory level and the number of customers in the pool. The problem is modelled as a semi-Markov decision problem and derived the optimal decision rule for the selection of pool customer rate using linear programming formulation.

	\subsection{\textbf{QIS with Retrial of Customers}} \label{subsec8}
	
	 Vijaya Laxmi and Soujanya  \cite{pv} investigated a  continuous review $(s, S)$ inventory system with negative customers and multiple working vacations at a service facility, in which a customer-demanded item is issued once the item has been serviced. When the inventory level becomes zero, the server takes multiple working vacations. Replenishment, vacation, and service times during regular busy periods and vacation periods are assumed to be exponentially distributed. The steady-state distribution of the model is computed using $MAM$. Assuming a suitable cost function for the inventory system, optimum values of $s$ and $S$ are obtained, minimising the total expected cost function.

	Kathiresan \textit{et al.} \cite{jkm} studied a retrial inventory system with negative Customers under continuous review $(s, Q)$ policy. The authors assumed that the customer chooses two types of service: ordinary and negative. The negative customer removes any of the ordinary waiting customers from the system, including the one at the service point. Each arriving customer may choose either a type  1 service or a type 2 service. A completion of type 1 service, one customer will receive service only, and in type 2 service, one customer receives service and reduces the inventory by one item. The steady-state behaviour and various system performance measures are derived, and the cost analysis is studied numerically. The model can be used to study different service types for the $(s, Q)$ ordering policy.

	A finite retrial inventory system with a service facility and multiple vacations for two heterogeneous servers who avail primary and retrial customers were considered by Suganya and Sivakumar  \cite{scs}. The demand process is according to $MAP$. Service durations for each server, lead time, and retrial time followed exponential distribution. The server becomes idle if there is no customer and the inventory level is positive. The idle period is terminated whenever the customers are from the orbit or primary. The stability of this system is analyzed, performance measures are computed, and the stationary distribution of the system state is obtained. The sensitivity of various system performance measures is studied using numerical illustrations.

	Jeganathan \textit{et al.}  \cite{jk} analyzed a Markovian $(s, Q)$ inventory-queueing system with server disruptions of two heterogeneous servers. The authors presumed the items involved in the service process to be non-perishable. The total capacity of the queue is assumed to be finite. If the queue is full, new customers arriving will enter into a finite-size orbit. When both the orbit and queue are full, the customer is lost. The matrix method determines a heterogeneous system's steady-state joint probability distribution of customers level in the queue, retrial group, server status, and stock level. A comparative study is done on the significant effect of a heterogeneous system with a homogeneous system and computed performance measures and appropriate cost functions. Also obtained that a heterogeneous system is more effective than a homogeneous system.

	Rejitha \textit{et al.}  \cite{krr} discussed two retrial inventory models having two modes of service rates under the $(s, Q)$ policy. The arrivals follow the Poisson process and will be entered into a buffer of finite capacity. Service time is supposed to follow an exponential distribution, and service is given at a reduced rate when the inventory level falls to $s$. When the buffer is full, an orbit of infinite capacity is provided to the customers, who can retry for service. The authors considered the lead time, inter-arrival, and inter-retrial times to be exponentially distributed and used $MAM$ to analyze the models. Various performance measures and a cost function are derived, and the two models are compared graphically.

	Jomy Punalal and Babu \cite{jpsb} focused on a single server retrial inventory system on self-generation of priorities with $(s, Q)$ order policy. A retrial inventory is considered and gave a probability rule for returning to orbit after an unsuccessful retrial. A customer moves to an orbit of infinite capacity when he finds the server is busy; if the server is free, a priority-generated customer gets immediate service. Service times are distributed exponentially, with differing rates for ordinary and priority consumers. The matrix analytic method is used to find the solution, a level-dependent quasi-birth-death process, and obtained important performance measures. An optimization problem is also discussed by introducing the profit function and calculating the optimum $(s, S)$ pairs.

	Jeganathan \textit{et al.} \cite{kjass} investigated the queue-dependent service rates in the stochastic QIS. If the server is accessible and there exists a positive stock in the stochastic QIS, an arriving customer immediately gets the service. Based on the classical retrial policy, the orbital customer can only compete for the service by enrolling in the waiting hall. The authors derived the generalized results for both homogeneous and non-homogeneous service rates. The service rate of any arrival depends on the number of customers in the queue. The Neuts Matrix Geometric approach derives the model's stability condition and stationary probability vector. The numerical examples explored the impact of queue-dependent service rates.

	A Markov model of a finite storage capacity QIS with primary, retrial, and feedback customers was considered by Melikov \textit{et al.}  \cite{maas}. One of the major assumptions in this paper is that the orbit for repeated customers could be formed by the feedback customers and the primary customers, with the destructive customers causing damage to items. Destructive customers do not require items as the inventory level instantly decreases by one upon their arrival. The authors adopted two replenishment policies: $(s, Q)$ and $(s, S)$, and developed a unified method to calculate the performance measures under the two replenishment policies.

	The stock-dependent Markovian demand of a retrial queueing system with a single server and multiple server vacations is studied by Sugapriya \textit{et al.} \cite{scnm} under a continuous review $(s, Q)$ ordering policy. The steady-state joint distribution of the number of customers in orbit, the server status, and the inventory level was obtained in the Matrix geometric approach with the Neuts-Rao truncation technique. The authors derived the significant system performance measures, and the waiting time distribution was investigated using the Laplace-Stieltjes transform and also discussed the merits of the proposed model, especially the important characteristics of the system with stock dependency over non-stock dependency.

	\subsection{\textbf{Production-inventory Systems}}\label{subsec9}
	
	In the production-inventory system (PIS), items can be replenished through the production unit. In most cases, the $(s, S)$ - policy is used in production-inventory models, wherein a production cycle starts when the inventory level drops to $s$ in response to a demand or deterioration. The production process stays in the "on-mode" until the inventory level hits the maximum level $S$, at which stage the production process is "switched off". 
 
 Yonit Barron \cite{ybar} considered a continuous-review PIS under the generalized order-up-to-level policy. The storage level was characterized as a reflected fluid process with uncertainties in the production and demand rates. The $(s, S)$ control policy is generalized to include order cancellations, exponential lead times, and unsatisfied demand. The cost associated with the model includes an order cost for each order or each cancellation, a purchase cost, a storage cost, and a penalty cost due to lost sales. The authors derived the explicit cost components of the resulting costs by taking a simple probability approach and applying stopping time theory to fluid processes and martingales. The cost components highlighted the impact of the system parameters on the different costs and were used for optimization purposes. Moreover, they obtained the optimal parameters using the closed-form expression of the components, and the system behaviour and its properties are studied and focused on the impact of the lost sales cost on the total cost.

	An $(s, S)$ production-inventory model with varying production rates and multiple server vacations was analyzed by Beena and Jose \cite{pbj}. The model considered customer arrival as $MAP$ and service time as $PH$ distribution. If there is at least one customer in the waiting area and the inventory level is positive, the service would occur at the end of a vacation period. They obtained several essential performance measures and defined a suitable cost function and an algorithmic solution for the problem using the $MAM$.

	A PIS under $(s, S)$ policy with unit items produced at a time, heterogeneous servers, vacationing servers, and retrial customers are discussed by Jose and Beena  \cite{kpjp}. The authors assumed that when the inventory level reaches zero or the orbit is empty, or both, one of the servers takes multiple, exponentially distributed vacation time. If the server finds the same at the end of the vacation, it immediately takes another vacation. During the stock-out period, server busy period, or vacation period, the demands that arrive enter an orbit of infinite size. When the customer in orbit cannot attempt due to the busy servers or the inventory level being  zero, he/she under Bernoulli trials may return to the orbit or leave the orbit. They obtained the ergodicity condition and used the $MAM$ to compute the steady-state probabilities. They computed several significant system performance measures.

	Jose and Reshmi \cite{jkprps} considered a continuous review perishable inventory system with a production unit and retrial facility, where customers arrive in a homogeneous Poisson distribution with $(s, S)$ ordering policy and perishable items. After every unsuccessful retrial, the customer returns to the orbit with a pre-allotted probability or is lost forever with complementary probability. Customers lost due to stock out can be considerably minimized because of the model's retrial facility. Finding the stationary probability vector using the matrix geometric method simplified the acquisition of some critical performance measures.

	Noblesse \textit{et al.}  \cite{amn} studied a continuous review finite capacity PIS with two products in inventory. The model reflects a supply chain that operates in an environment with high levels of volatility. The authors considered inventory replenishment as an independent order-up-to $(s, S)$ policy or a can-order $(s, c, S)$ joint replenishment policy in which the endogenously determined lead times drive the parameters of the replenishment policy. The production facility is modelled as a multi-type $MMAP[K]/PH[K]/1$ queue. The steady-state distribution is derived using the $MAM$. The inventory parameters that minimize the total ordering and inventory-related costs are obtained using numerical methods.

	Under a base stock strategy for inventory control, a single-item PIS is investigated by Cruz \textit{et al.} \cite{dch}. The authors used a Gordon Newell network with $z \geq1$ customers and $K + 1$ exponential single server stations $J = \left\{1, \dots, K, K + 1\right\}$ to represent the complete PIS. The (replenishment) production system is represented by the stations $j = 1, \dots, K $. In contrast, the finished goods inventory is represented by the station $K + 1$. The base stock level, the total of the finished goods and work-in-process inventory, determines the network's population size. With each exogenous demand that follows a Poisson process, a production order for a new unit is released, and the amount of the work-in-process inventory rises. They deal with the issue of choosing a base stock level that minimizes a standard cost function. This cost function includes the holding costs for inventory items, the production costs for replenishment system items, and the costs for lost sales. Even though the system is transformed into a standard closed Gordon- Newell network, it is easy to understand with respect to steady-state distribution and expected queue lengths (first-order properties). The study revealed that proving the second-order properties of interest is quite challenging.

	Otten and Daduna \cite{od} analyzed a PIS with two classes of customers with different priorities admitted into the system via a flexible admission control scheme. The service time is exponentially distributed with parameter $ \mu> 0$ for both types of customers. An arriving demand that finds the inventory depleted is lost because the inventory management follows a base stock policy (lost sales). To know the equilibrium behavior of the system, the global balance equations of the related Markov process are examined and deduced structural features of the steady-state distributions are. Using the Foster-Lyapunov stability criterion, they derived a sufficient condition for ergodicity for both customer classes in the case of unbounded queues. If the customer is at the head of the line and the inventory is not depleted, the service begins immediately. Otherwise, the service starts instantaneously when the next replenishment arrives in the inventory. The queue admission control is flexible and includes two parameters $(s; p)$ that primarily limit the entry of low-priority customers. Although the authors proved that the condition is sufficient and necessary for the parameter constellation $(s; 1)$, there is no assurance for the general case.

    Kocer and Ozkar \cite{ko} examined a QIS where the server is subject to breakdown with two types of customers, Type 1 and Type 2 are priority and ordinary customers, respectively. Type 1 customers have a non-preemptive priority over Type 2 customers. A continuous review $(s, S)$ inventory policy has been established, and inventory is supplied through inner production. Due to breakdowns, service may be interrupted at any time. Minor or major breakdowns may occur, and following the recovery process, the priority customer class receives the service first. A major breakdown occurs with probability $(1-p)$, and a minor breakdown with probability $(p)$. Only when a server is busy might it be interrupted due to a breakdown. A possible interruption does not affect an idle server. Because of the service interruption, no inventory is lost. Customers who are being served at the time of a service interruption queue up until the server is repaired. The time required to repair a breakdown depends on its nature. They modelled the system as a level-dependent QBD model and computed steady-state probabilities and an optimal strategy is investigated numerically.

    \subsection{\textbf{QIS with Random Environment}}\label{subsec10}
The behavior of the system is determined by the current state of the Markovian random process with a finite state space. A "random environment" is a term that describes the method of Markovian random process under a fixed state.

    Anbazhagan \textit{et al.} \cite{aavaj} proposed a multi-server queueing-inventory system in a random environment. The state of the random environment determines the behavior of the system. The demands arise from a random environment with $R$ various modes, each requiring a single item unit, and the system has a maximum inventory capacity of S units. The customer arrival follows the Markovian arrival process, the service time follows a phase-type distribution, and the inventory replenishment is instantaneous. The parameters of the arrival and service flow change whenever the states of a random environment change. The authors obtained the joint probability distribution of the inventory level and the number of customers in the steady-state case. They derived the sojourn time distributions of arbitrary customers using the Laplace-Stieltjes transform.

    Jacob \textit{et al.} \cite{jjsdk} introduced the concept of optional items for service in a finite number of randomly changing environments, with a single server multi-commodity QIS with positive service time. The authors assumed that there is only one essential item under the $(s, S)$ ordering policy, and all other items are to be optional. The authors defined a new concept called Marked Markovian Environment Change Process (MMEAP). The authors further assumed that the customer's demand follows $MAP$, thus, the environmental change is turned to be an MMEAP, the service time is $PH$ distributed for essential inventory, and those optional items are distributed exponentially, and they analyzed the stability of the system. A cost function  is constructed and investigated optimal control policies for the different types of commodities.

    Jacob and Krishnamoorthy \cite{jjk} investigated a single commodity queueing- inventory with $n$ different environments. The customer's arrival follows the Poisson process. Service time and lead times are mutually independent, exponentially distributed random variables. The sojourn time follows an exponential distribution. No research has been published on product form solutions of queueing inventory systems that evolve through a finite number of random environments. The authors obtained asymptotic product form solutions in 2 variants by imposing certain restrictions. For the first case, the customer arrival rate $\lambda_{i}=\lambda$, service completion rate $\mu_{i}=\mu $ and the replenishment rate $ \beta_{i} =\beta$, and the control policy $(s_{i}, S) = (s, S)$, $1 \leq i \leq n$ are all assumed to be the same for all n different environments and the transition rate from $(i, j, m_{1}) \rightarrow (i, j, m_{2})$ that is due to the change in environment by retaining the same inventory level is fixed to be $\theta_{m_{1}} P_{m_{2}}$ with these assumptions the authors obtained the asymptotic product form solution as 
\begin{eqnarray*}
 \mbox{\boldmath{$\Tilde{x}$}}(i, j, k) = \rho^{i} (1 - \rho) \Tilde{\zeta} (k),\ i \geq 0,
\end{eqnarray*} where  $\mbox{\boldmath{$\Tilde{x}$}} = ({\mbox{\boldmath{$\Tilde{x}$}}_{0}}, {\mbox{\boldmath{$\Tilde{x}$}}_{1}}, {\mbox{\boldmath{$\Tilde{x}$}}_{2}}, \dots).$
For the second variant, the authors assumed that only cyclic environment change is allowed and the transition rate from $(i, j, m_{1}) \rightarrow (i, j, m_{2})$ that is due to the change in environment by retaining the same inventory level is fixed to be $\theta_{m_{1}}$, the asymptotic product form solution is,
\begin{eqnarray*}
 \mbox{\boldmath{$\hat{x}$}}(i, j, k) = \rho^{i} (1 - \rho) \hat{\zeta}(k), \ i\geq 0,
\end{eqnarray*} where $\mbox{\boldmath{$\hat{x}$}} = ({\mbox{\boldmath{$\hat{x}$}}_{0}}, {\mbox{\boldmath{$\hat{x}$}}_{1}}, {\mbox{\boldmath{$\hat{x}$}}_{2}}, \dots)$.

    Vinitha \textit{et al.} \cite{vaaj} discussed the random environment with two classes of suppliers and impulse customers. The authors considered two categories of suppliers, temporary and regular suppliers, and both lead times follow an exponential distribution. These two types of suppliers have a significant impact on inventory management and help to reduce replenishment difficulties. The arrival of customers follows the Markovian arrival process from $N$ different states of the random environment. They order $Q_{1}$ items from the temporary supplier whenever the inventory level approaches $r$, where $Q_{1} = S -r$. The parameters of the arrival flow and retrying flow are simultaneously changed if the states of a random environment are changed. Customers who only require one item can select any arrival mode (state) RE. The status of the RE determines the behavior of the model. When the inventory level is positive, the service process of the system is assumed to be instantaneous. The authors obtained the steady-state probability distributions using $MAM$. Finally, the authors proved that the value of the setup cost and the lead time rate are both small, the best supplier is a temporary supplier, and both have higher value, making the best supplier a regular supplier.

    Sonja Otten \textit{et al.} \cite{srhk} considered exponential single server queues with state-dependent arrival and service rates that change in response to external environmental factors. The authors assumed a bidirectional interaction, that is,  the state of the environment influences queue transitions, while the movements of the environment depend on the status of the queues and obtained product form steady-state distributions. The authors developed ergodicity and exponential ergodicity criteria using Lyapunov functions. They presented how to construct two-sided bounds for the throughput using the example of a coupled production-inventory system for nonseparable systems. The queue and the environment processes are decoupled asymptotically and in a steady state. Using examples, the authors demonstrated the principles for setting upper and lower bounds on the departure rates of customers serviced (throughputs) of nonseparable systems based on throughputs of related separable systems. Another path of investigation into nonseparable queueing-environment systems is appropriately changing the system's transition rate matrix (possibly only the environment coordinate) to get a product form approximation of the stationary distribution.

    \subsection{\textbf{QIS - Game-theoretic Analysis}}\label{subsec11}
    Until the early nineties of the last century, researchers were focused on investigating the real-life problems of administrative measures for the optimization of queues. More precisely, they assume arriving customers would always join a queueing system even if their expected waiting time is high. Nevertheless, this is not always realistic. In this case, a more realistic assumption is that the customers are willing to receive quality service, but at the same time, they are not willing to wait for a long time in the queue. Thus, customers are free to make decisions about their actions in the system with the objective of maximizing their own benefit or social benefit. The study of queueing systems with customers having the option to make decisions was first investigated by Naor \cite{naor} in 1969. After Naor’s pioneer work, there is an emerging tendency to investigate the behaviour of customers in queues who have the freedom to make decisions taken place. For a detailed survey of strategic queueing systems, refer to Hassin and Haviv \cite{hassin03}, and Stidham Jr. \cite{stidham}. However, the study of social/individual optimal strategy in queueing-inventory systems is not well discussed. In this section, we discuss the articles on the game-theoretic analysis of QIS with MTS and MTO policies, we refer to the monograph by Rafael Hassin \cite{Rhassin} for a partial discussion on the same.

	Li \textit{et al.}  \cite{Guo} considered an $M/M/1$ make-to-stock, finite-capacity production system with high setup costs and delay-sensitive customers. To the best of our knowledge, this is the first reported work on QIS from a game-theoretic perspective, in which the model is well integrated with three different research areas queues, inventory, and game theory. The production manager uses a two-critical-number control policy to balance setup and inventory-related expenses. According to this policy, production begins when $N$ customer orders accumulate and stop immediately when the inventory level reaches $S$, where $S\geq0$. The customers decide whether to stay for the product or leave without purchasing based on utility values determined by the production manager's control decisions. They frame the problem as a Stackelberg game involving the production manager and the customers, with the former serving as the game leader. In the former case, there may be two positive equilibria, with only the larger one remaining stable. In the latter case, there is only one stable positive equilibrium effective arrival rate. They also demonstrated that reducing the critical number of waiting for customers or increasing the critical number of cumulated inventory can reduce the expected waiting time, with the former having a greater marginal impact than the latter. The numerical results reveal that there may be lost sales when demand traffic is very light or heavy, and the optimal cost rate can be significantly lower when dealing with non-strategic customers.

	Jinting Wang \textit{et al.}  \cite{Wang1} investigated the customer's individual optimal and socially optimal strategies of an $M/M/1$ QIS and analyzed the optimal pricing issue that maximizes the server's revenue. This is the first contribution of interest for analyzing the equilibrium behaviour of customers along with the pricing mechanism in an $M/M/1$ service-inventory system. Here the inventory status is a system parameter taking value from $0$ to $Q$ ($Q$ is the capacity of the inventory), which is controlled by a predetermined continuous review policy. The authors obtained individual and socially optimal strategies in fully unobservable and partially observable cases. By introducing the optimal price that maximizes the server's revenue, the maximum revenue in a partially observable case is not always more significant than that in a fully unobservable case. They numerically demonstrate that the revenue maximum is equal to the social optimum in the fully unobservable case in most cases.

	Wang \textit{et al.} \cite{Wang2} is the first to combine a two-stage inventory model with a loss-averse corporation with a reference point and strategic customers with declining valuation. This paper considered a two-stage decision model and examined the effect of the reimbursement contract and price commitment strategy on the alleviation of strategic customer behaviour. Consumers can choose whether to buy at full price in Stage 1 or to wait for the salvage price in Stage 2. They might not be able to purchase the product if it runs out in stage 2. The firm's goal is to select a base stock policy and find an optimal price to maximize its expected utility. Customers' goal is to decide whether to buy or wait strategically to maximize their payoffs. The authors formulate the problem as a Stackelberg game between the firm and the strategic customers, where the leader is the firm. They obtained that the reimbursement contract cannot alleviate the negative influence of strategic customer behaviour, and the firm's profits may increase by providing price commitment strategies. Also, it is analyzed numerically.

	Zhang  \& Wang  \cite{Zhang} considered a make-to-stock production system with one product type, dynamic service policy, and delay-sensitive customers. Corresponding the number of customers is zero or not, the system operates at a low or high production rate. The model is formulated as a Stackelberg game between the production manager and customers, the former is the game leader, and the leader seeks an appropriate threshold to obtain the maximum revenue. This work focused on two kinds of strategies: the customer's equilibrium strategy, the Nash symmetric equilibrium strategy, and the manager's optimal strategy. Numerical results show that the expected cost function in an observable case is not greater than that in an unobservable case. These two cases are identical if a customer's delay sensitivity is relatively small. A demarcation line is depicted to differentiate the MTS policy from the MTO policy, and the optimal inventory threshold is found within the MTS region.

	Cai \textit{et al}  \cite{Cai} examined the pricing and inventory choices in a make-to-stock system from the view of a profit-maximizing server and a social planner with delay-sensitive customers. Almost unobservable and fully observable information scenarios are the main focus of this work, in which customers decide whether to purchase the item based on its utility. Profit maximizing and welfare maximization were the two approaches used to study the two-dimensional decision problem in this work, which resulted in different pricing and inventory strategies in both observable and unobservable queues. They solved the profit maximization problem by formulating it as a Stackelberg game between the server and customers. The authors obtained the numerical results that the optimal price is sensitive to delay information, whereas the optimal stock level is not.

	Lee \textit{et al.}  \cite{Lee} explores a market with a duopoly between two firms, each of whom runs its service inventory systems. The two retailers sell different brands of similar products. Both competitors set their prices to maximize profit while accounting for inventory holding costs, ordering costs, and the cost of lost sales. Each competitor's service inventory system is modelled as an $M/M/1$ queue with varying arrival rates that depend on price. Three pricing games were developed and examined using a game theory approach: a parallel pricing game, a sequential pricing game, and a unified pricing game. Two competing retailers set their retail prices simultaneously in the parallel pricing game. The players in this strategic game are the retailers, who are fully informed. From this pricing game, they obtained that a retailer's price increase reduces the arrival rate to their service-inventory system, and the replenishment cycle lengthens as the arrival rate decreases. The sequential pricing game is modelled after the duopoly Stackelberg game, in which one retailer is the pricing leader, and the other is the pricing follower. In the uniform pricing game, it is assumed that a monopolist sells substitutable products and that they are jointly determined to maximize profits. Several numerical examples reviewed the impact of modifications in various parameters on each retailer's pricing game.

	Kim \textit{et al}. \cite{Kim} considered a production–inventory queueing system, arriving customers can observe the system state and then decide whether to wait or not to wait for the product. In the beginning, they demonstrated a customer equilibrium strategy for a certain joint price and inventory control. The production manager can determine joint pricing and inventory control. Finally, they showed the existence of the customers’ equilibrium strategy and studied when it is unique. The equilibrium is unique if the production time distribution decreases the mean residual life. They also looked into profit and social welfare maximization and discovered that the optimal joint pricing and inventory control policy maximized the social benefit rate.

	\subsection{\textbf{QIS with Different Scenarios}}\label{subsec12}
	
	There are a few more articles on QISs with various cases by several researchers that do not fall under the above-discussed categories. We feel that, without mentioning those important contributions to QISs, this survey will be incomplete. This section is mainly dedicated to providing a detailed literature review on the QISs with various cases.

 Al Maqbali \textit{et al.} \cite{maq} introduced a novel approach involving the batch arrival of customers at a transport station. They introduced the concept of designated capacity reserved vessels for different customer categories, a previously unexplored idea. The authors applied the Batch Marked Markovian Arrival process (BMMAP) to treat each customer category as a unified entity. This system encompasses two types of rooms: waiting rooms at the transfer station and service rooms on the vessel, specifically catering to customers of type J. In this model, an Erlang clock governs the phase for initiating customer searches for the current vehicle. When the transport ship arrives at the station, an Erlang clock with order m begins counting. This clock triggers the scheduling of the next vehicle when it reaches stage L. This scheduling entails an exponentially distributed time interval required for the vehicle to reach the station. The vehicle was already present in the station leaves the moment given by min{realization of the Erlang clock, the arrival of the next vehicle, all seats are filled}. They derived the stability condition.

	Melikov \textit{et al.}  proposed a new restocking policy for QISs in \cite{am}. They used the fact that order sizes are variables that depend on current stock levels as a distinguishing aspect of this policy. The authors considered two models: one with instantaneous service and the other with nonzero service time. The explicit formulas for calculating the characteristics of QIS models are obtained with instantaneous service of spending requests when utilizing the proposed policy, variable size order policy (VSO), and the well-known $(s, S)$ policy of two levels. The authors assumed that impatient spending requests can form a finite or infinite length queue in models with nonzero service time and that the intensity of spending request loss from the queue depends on the current stock level. In the case of utilizing the proposed restocking policy, they also established accurate and approximate methods for determining the characteristics of the QISs under consideration. They compared the properties of the QISs under consideration in the case of various restocking policies.

	G Hanukov \textit{et al.}  \cite{gh} investigated an innovative approach to increase the efficiency of queueing systems by utilizing the server's idle time to produce preliminary services for future incoming customers. They used a single server Markovian queue and constructed a two-dimensional state space that considered queue sizes and inventory levels.
	Assuming linear costs for customers waiting in line and for stored preliminary services, a cost analysis determined the optimal maximal number of stored preliminary services in the system. Using probabilistic methods, the  authors obtained the steady-state probabilities of the system states and various performance measures. The advantages of the approach, in terms of cost savings, as compared with the classical $M/M/1$ model, are illustrated using numerical examples with graphs.

	Varghese and Shajin \cite{vd} considered an $M/M/1/S$ inventory model with a finite storage system under $(s, S)$ policy. The service time and lead time follow the exponential distribution, and the arrival follows a non-homogeneous Poisson process. They assumed that a new arrival joins the queue only if the inventory is not less than the number of demands in the system. The system characteristics are evaluated using the $MAM$. The authors constructed a cost function to numerically investigate the optimal pair $(s, S)$.

	Boxma \textit{et al.}  \cite{bojp} studied the relationship between the $(S-1, S)$ inventory model and three well-known queueing models: the Erlang loss system, the machine repair model, and a two-node Jackson network. This paper's main objectives were to (i) highlight the connections between the $(s-1, s)$ inventory model and various important queueing models, particularly the connections between the $(s-1, s)$ model and the Erlang loss system, and (ii) use the latter connections to obtain results for the $(s-1, s)$ model's key performance measures, such as the so-called virtual outdating time and the number of items on the shelf. Investigating whether one could once more take advantage of the relationship to queueing models while taking into account variations and generalizations of $(s-1, s)$ is interesting. This relationship made it possible to collect performance measures like the distribution of the empty shelf duration and the density of the virtual out-dating time.

	Choi \textit{et al.}  \cite{ckhy} discussed a QISs with $PH$ service distributions and discusses their contributions to inventory management. The authors classified the models in the literature based on the queueing model and features such as replenishment policy, stock-out assumption, shelf time, and positive service time. This review provides a clear overview as a starting point for further study of a queueing inventory model.

	The impact of a sequential supply system on the total cost of a two-echelon inventory system with one central warehouse, some retailers, lost sales and emergency shipment is studied by Tayebi \cite{th}. They assumed that if stock is out at the warehouse, the external supplier replenishes retailers' stocks by emergency shipments. The warehouse and retailers applied a base-stock ( $(S-1, S)$) ordering policy, and the unsatisfied demands at retailers were considered lost. The lead times of retailers are assumed to be Erlangian distributed. Under the realistic assumption that the replenishment orders do not cross in time, the Erlangian lead times become stochastically dependent, which affects the optimal base-stock levels at retailers and, consequently, the optimal base-stock level at the warehouse. The authors estimated the warehouse's demand rate, computed the inventory system's total cost function and presented an algorithm to compute the optimal base-stock levels at the warehouse and retailers.
	
	Two QISs with batch $MAP$ demands, $PH$-service times, and $(s, S)$ replenishment policy were studied by Chakravarthy and Rumyantsev \cite{src}. The authors assumed that an arriving customer is lost if the inventory level is zero in one model, whereas in the second model, a customer can be lost either at arrival or at the time of service completion. The second model removes all waiting customers from the system due to zero inventory. They employed the classical $MAM$'s for one model, and in the case of the second model, they incorporated simulation using ARENA software. The authors presented an optimization problem to compare both models and investigate the nature of optimum configurations for the two models separately.
	
	Keerthana \textit{et al.}  \cite{kms} analyses a stochastic inventory system with a service facility in which only one customer arrives at a time, intending to find the optimal inventory level and reorder point that minimizes the long-run expected cost rate. They assumed the customer arrived according to a renewal process and demanded the item delivered to the customer after performing an exponentially distributed service time. 
	They adopted the $(s, S)$ ordering policy with exponentially distributed lead times. They derived the stationary probability distribution for the number of customers in the system, inventory level at arrival epoch, and arbitrary time points. The authors provided a sensitivity study to illustrate the effects of parameters and cost on the optimal values.
	
	Chakravarthy and Hayat  \cite{csh} introduced the concept of multiple vendors responsible for replenishing inventory so that replenishment occurs using two vendors as opposed to the classical approach of using only one vendor. Under the assumptions of a two-vendor system wherein the demands occur according to a $MAP$, the service times are $PH$, and the lead times are exponentially distributed with a parameter depending on the vendor, the authors analyzed the model in steady-state using the $MAM$. Interesting numerical examples are presented, including comparing the one and two-vendor systems.

	Y Zhang \textit{et al.}  \cite{yzdy} considered a QIS model with the server’s multiple vacations and impatient customers, in which the new customer arrives and waits for service when the server is off due to vacation. Using the method of truncated approximation, they obtained the matrix geometric solution of the steady-state probability and derived some relevant performance measures. The effect of the probability and impatient rate on some performance measures was investigated using numerical analysis. The authors computed the optimal policy and cost using the genetic algorithm and derived the optimal service rate.

	Chakravarthy and Rao \cite{csrrb} proposed an opportunistic-type inventory replenishment model with two models: Opportunistic Model 1 and Opportunistic Model 2. Customer arrivals follow the $BMAP$. In Model 1, the customers are lost when the inventory level is zero (irrespective of whether the service is idle or busy). In Model 2, the customers are allowed to enter even when the inventory level is zero, but as long as the server is busy serving. This article thoroughly examines system behavior to compare the conventional inventory system with the opportunistic models and finds the optimum strategy to manage inventories in a more broad situation using $MAP$ for the demand process and phase type distributions for processing.

    Samanta \textit{et al.}  \cite{siv} examined a continuous review $(s, Q)$ inventory system with a service facility. They consider an $M/G/1$ QIS with a finite waiting capacity $N$. The demands arrived according to the Poisson process. Every customer requires a single product with an arbitrary and unexpected service period. A new customer will be accepted or rejected based on the number of customers enrolled in the system. Using the embedded Markov chain technique, the authors computed the joint probability distribution of the number of customers in the system and the number of items in inventory during the post-departure period. The authors developed an analytical expression for the average total cost of the operating system. In addition, the authors discussed a few simple numerical examples showing how the appropriate order quantity depends on the utilization factor, service rate, and arrival rate. They revealed that the overall cost varies more quickly when the utilization factor exceeds one.

    Sebastian \textit{et al.}  \cite{sass} presents a detailed analysis of a finite-source inventory system with a service facility and postponed demands. The inventory system with a service facility that consists of a single server and a finite waiting hall of size $K$, which includes the customer who receives the service. The arrival of customers for unit items follows a quasi-random process. The service time to process the item follows phase-type distribution and replenishment following the $(s, S)$ order policy. According to the Bernoulli trial, a newly arriving customer who finds the waiting hall occupied can either join the pool for postponed customers or leave the system. The joint probability distribution of the inventory level, the number of customers in the pool, and the number of customers in the waiting hall are obtained in the steady-state case. The authors presented a steady-state analysis of the inventory model under the assumptions of quasi-random arrival processes, phase-type service times, phase-type lead times, and prefixed selection rules for the selection of postponed customers in the pool. The SCV (squared coefficient of variation) of the service time and lead time distributions influences the optimal values. They concluded that SCV is essential for deciding ordering policies, selection policies, and waiting room size as a result of the numerical studies. As a result, they minimize the cost or waiting time associated with lead times or service times by controlling the SCV.

    Mohammad Rahiminia \textit{et al.} \cite{Mrahim} are the first to design a mathematical model of healthcare and waste management using a data-driven method. Using a bulk service and two different multiserver queueing models, the authors present a non-linear multi-shift mathematical model to control the congestion in the medical centre and medical waste during the pandemic. After categorising patients based on their current health condition, they used $M/M/C$ queueing systems to control patient congestion. They then designed each class as an $M/M/C$ queueing system to simulate real-world pandemic conditions. This is the first time the HWQIS connected with the medical centre has used a single server Markovian queue with bulk service $M/M[y]/1$. In this approach, they combine a balk queuing system for inventory waste generated by the medical centre with two queuing systems for patients in the medical centre.

    \section{Conclusion and Future Work} 
    This survey provides a comprehensive overview of the contributions of queueing-inventory systems from 2016 to 2023 by highlighting various applications in real-life problems. This survey also comprises the contributions of the gradual developments of the product-form solutions of the QIS, the current trends of the research topics, such as game-theoretic analysis of QIS and QIS with the random environment, which has not been previously discussed in the earlier survey articles on QIS (See \cite{klm}, \cite{ksn} \cite{ckh}). Furthermore, discrete-time analysis of QIS has received less attention; future works include analysis of discrete-time QIS, stochastic control problems in dynamic QIS, the interdependence of QIS through a semi-Markov process, and economic analysis of QIS through a game-theoretic approach can open up further research opportunities. However, the challenges regarding the transient solutions of queueing-inventory systems remain a gap in the literature because of the model's complexity.

    \begin{landscape}
			\begin{table}\vspace{1cm}
		\begin{tabular}{|p{1.8cm}||p{0.9cm}|p{0.6cm}|p{0.9cm}|p{0.7cm}|p{0.4cm}|p{0.7cm}|p{0.6cm}|p{0.8cm}|p{0.7cm}|p{0.7cm}|p{0.75cm}|p{0.7cm}|p{0.65cm}|p{0.7cm}|p{0.7cm}|p{0.6cm}|p{0.8cm}|p{0.7cm}|p{0.8cm}|}
			\hline
			\multicolumn{20}{|c|}{Keywords} \\
			\hline
			Author&	Poisson Arrival& MAP&		MMAP& Expo- nential Servi- ce Time&	PH Service& Prod- uct Form Solu- tion&	Peri- shable Inve-ntory&	Multi- comm odity&	Post poned demand&	Impat- ient Customer&	Single Server&	Multi Server&	Retr-ial Que-ues&	MAM&Conti- nuous Review&	Discr- ete QIS&	Produ- ction Inve- ntory&Server Vaca- tion&	Multi-class cust- omers \\
			\hline
			Khalid A Alnowibet \cite{akaa}   & $\checkmark$    &     & &$\checkmark$ & &$\checkmark$ & & & &$\checkmark$ &$\checkmark$ & & & & & & & & \\
			\hline
			M Amirthakodi \cite{ams}& $\checkmark$ &  & &$\checkmark$ & & & & & & &$\checkmark$ & &$\checkmark$ & &$\checkmark$ & & & & \\
            \hline
			N Anbazhagan \cite{aavaj}&&$\checkmark$   & & &$\checkmark$  & & & & & &&$\checkmark$  && && & & & \\
			\hline
			Anilkumar \cite{amp}    & &  & & & & & & & & &$\checkmark$ & & &$\checkmark$ & &$\checkmark$ & & &  \\
			\hline
			Ann M \cite{amn} &  &  & & & & & &$\checkmark$ & & & & & &$\checkmark$ &$\checkmark$ & &$\checkmark$ & & \\
			\hline
			A Z Melikov \cite{azmm} & &  & &$\checkmark$ & & &$\checkmark$ & & & &$\checkmark$ & & & & & & & & \\
			\hline
			J Baek\cite{jb}&  &  &$\checkmark$ & & & & & & &$\checkmark$ &$\checkmark$ & & & & & & & &$\checkmark$ \\
			\hline
			Y Barron \cite{bar} &&&&$\checkmark$&&&$\checkmark$&&&&&&&&$\checkmark$&&&& \\
			\hline
			Y Barron \cite{by} &&&&&&&$\checkmark$&&&&&&&&$\checkmark$&&&& \\
			\hline
			O Baron \cite{bp} &$\checkmark$&&&$\checkmark$&&&$\checkmark$&&&&&&&&$\checkmark$&&&& \\
			\hline
			B Benny \cite{bbv} &$\checkmark$&&&$\checkmark$&&&&&&&&$\checkmark$&&&&&&&$\checkmark$ \\
			\hline
			M Bhuvanes-hwari \cite{bm} &&&&$\checkmark$&&&$\checkmark$&&&&&&&&$\checkmark$&&&& \\
			\hline
			O J Boxma\cite{bojp} &$\checkmark$&&&&&&&&&&&$\checkmark$&&&&&&& \\
			\hline
		\end{tabular}
	\caption{QIS with different categories}
	\end{table}
    \end{landscape}

    \small
	\begin{landscape}
		\begin{table}
		\begin{tabular}{|p{2cm}||p{0.9cm}|p{0.5cm}|p{0.8cm}|p{0.9cm}|p{0.4cm}|p{0.7cm}|p{0.7cm}|p{0.8cm}|p{0.75cm}|p{0.8cm}|p{0.8cm}|p{0.9cm}|p{0.9cm}|p{0.6cm}|p{0.6cm}|p{0.7cm}|p{0.9cm}|p{0.77cm}|p{0.9cm}|}
			\hline
			\multicolumn{19}{|c|}{Keywords} \\
			\hline
			Author&	Poisson Arrival& MAP&	BMAP& Expon- ential Servi- ce Time&	PH Service& Prod- uct Form Solu- tion&	Peris- hable 
			Inven- tory&	Multi- comm odity&	Impat- ient Customer&	Single Server&	Multi Server&	Retrial Queues&	MAM&	Conti- nuous Review&	Discr- ete QIS&	Produ- ction Inve- ntory&	Server Vaca- tion&	Multi-class cust- omers\\
			\hline
			S R Chakravarthy \cite{chasr} &$\checkmark$&&&$\checkmark$&&&&&&$\checkmark$&&&$\checkmark$&&&&& \\
			\hline
			Srinivas R Chakravarthy \cite{csh} &&$\checkmark$&&&$\checkmark$&&&&&&&&$\checkmark$&&&&& \\
			\hline
			S R Chakravarthy \cite{csr} &$\checkmark$&&&$\checkmark$&&&&&&$\checkmark$&&&$\checkmark$&&&&& \\
			\hline
			S R Chakravarthy \cite{csrrb} &&&$\checkmark$&&$\checkmark$&&&&&&&&$\checkmark$&&&&& \\
			\hline
			K H Choi\cite{ckhy} &&&&&$\checkmark$&&&&&&&&&&&&& \\
			\hline
			NN de la Cruz \cite{dch} &$\checkmark$&&&&&&&&&&&&&&&$\checkmark$&& \\
			\hline
                H Daduna \cite{dad} &$\checkmark$&&&$\checkmark$&&$\checkmark$&&&&$\checkmark$&&&&&&&& \\
			\hline
			A Dudin \cite{dakv} &&$\checkmark$&&&&&&&&$\checkmark$&&&&&&&& \\
			\hline
			M Fathi\cite{fmkm} &$\checkmark$&&&&&&&&&&&&&&&&&$\checkmark$ \\
			\hline
			Gabi Hanukov\cite{gh} &$\checkmark$&&&$\checkmark$&&&&&&$\checkmark$&&&&&&&& \\
			\hline
			Gabi Hanukov\cite{gb} &$\checkmark$&&&$\checkmark$&&&$\checkmark$&&&$\checkmark$&&&&&&&& \\
			\hline
			G Hanukov \cite{Hanv} &$\checkmark$&&&$\checkmark$&&&&&&&$\checkmark$&&&&&&& \\
			\hline
			Padmavathi \cite{ip} &&$\checkmark$&&&$\checkmark$&&&&&&&&&&$\checkmark$&&$\checkmark$& \\
            \hline
			J Jacob\cite{jjk} &$\checkmark$&&&$\checkmark$&&&&&&$\checkmark$&&&&&&&& \\
			\hline
			J Jacob\cite{jjsdk} &&&&&&&&$\checkmark$&&$\checkmark$&&&&&&&& \\
			\hline
			K Jeganathan \cite{jk} &&&&&&&&&&&$\checkmark$&$\checkmark$&$\checkmark$&&&&$\checkmark$& \\
			\hline
			K Jeganathan \cite{jkma} &$\checkmark$&&&$\checkmark$&&&$\checkmark$&&&&$\checkmark$&&$\checkmark$&$\checkmark$&&&$\checkmark$& \\
			\hline
			K Jeganathan \cite{jkrma} &$\checkmark$&&&$\checkmark$&&&&$\checkmark$&&$\checkmark$&&$\checkmark$&&&&&& \\
			\hline	
		\end{tabular}	
	\caption{QIS with different categories}
 \end{table}
    \end{landscape}

    \begin{landscape}
		\begin{table}
		\begin{tabular}{|p{1.8cm}||p{0.9cm}|p{0.55cm}|p{0.8cm}|p{0.9cm}|p{0.4cm}|p{0.7cm}|p{0.7cm}|p{0.7cm}|p{0.7cm}|p{0.7cm}|p{0.7cm}|p{0.7cm}|p{0.6cm}|p{0.65cm}|p{0.65cm}|p{0.65cm}|p{0.85cm}|p{0.77cm}|p{0.8cm}|}
			\hline
			\multicolumn{20}{|c|}{Keywords} \\
			\hline
			Author&	Poisson Arrival& MAP&	BMAP&	 Expon- ential Servi- ce Time&	PH Service& Prod- uct Form Solu- tion&	Peris- hable 
			Inven- tory&	Multi- comm odity&	Post poned demand&	Impat- ient Customer&	Single Server&	Multi Server&	Retr-ial Que-ues&	MAM&	Cont- inuous Review&	Discr- ete QIS&	Produ- ction Inve- ntory&	Server Vaca- tion&	Multi-class cust- omers\\
			\hline
			K Jeganathan \cite{jkvsh} &&&&&&&$\checkmark$&&&&$\checkmark$&&&&&&&&$\checkmark$ \\
			\hline
			J S A Jenifer \cite{jsa} &$\checkmark$&&&$\checkmark$&&&&&$\checkmark$&&$\checkmark$&&&&$\checkmark$&&&& \\
			\hline
			J Kathiresan\cite{jkk} &$\checkmark$&&&$\checkmark$&&&&&&$\checkmark$&&$\checkmark$&$\checkmark$&&$\checkmark$&&&$\checkmark$& \\
			\hline
			J Kathiresan\cite{jkm} &$\checkmark$&&&$\checkmark$&&&&&&&&&$\checkmark$&&$\checkmark$&&&&$\checkmark$ \\
			\hline
			Jomy Punalal\cite{jpsb} &$\checkmark$&&&$\checkmark$&&&&&&&$\checkmark$&&$\checkmark$&$\checkmark$&&&&&$\checkmark$ \\
			\hline
			K P Jose\cite{jkprps} &$\checkmark$&&&$\checkmark$&&&$\checkmark$&&&&$\checkmark$&&$\checkmark$&$\checkmark$&&&$\checkmark$&& \\
			\hline
			J W Baek\cite{jwbaek} &$\checkmark$&&&$\checkmark$&&$\checkmark$&&&&&$\checkmark$&&&&&&&& \\
			\hline
			M Keerthana \cite{kms} &&&&$\checkmark$&&&&&&&$\checkmark$&&&&&&&& \\
			\hline
			K Jeganathan\cite{kjm} &$\checkmark$&&&$\checkmark$&&&$\checkmark$&&&&&$\checkmark$&&&&&&$\checkmark$&$\checkmark$ \\
			\hline
			K Jeganathan\cite{kjna} &$\checkmark$&&&$\checkmark$&&&&&&&&$\checkmark$&&&&&&$\checkmark$&$\checkmark$ \\
			\hline
			K Jeganathan\cite{kjass} &$\checkmark$&&&$\checkmark$&&&$\checkmark$&&&&$\checkmark$&&$\checkmark$&&&&&& \\
			\hline
            U U Kocer\cite{ko} &$\checkmark$&&&$\checkmark$&&&&&&&$\checkmark$&&&&$\checkmark$&&$\checkmark$&&$\checkmark$ \\
			\hline
			
		\end{tabular}
	\caption{QIS with different categories}
	\end{table}
    \end{landscape}

    \begin{landscape}
			\begin{table}
		\begin{tabular}{|p{2cm}||p{0.9cm}|p{0.55cm}|p{0.8cm}|p{0.9cm}|p{1cm}|p{0.4cm}|p{0.7cm}|p{0.7cm}|p{0.8cm}|p{0.7cm}|p{0.7cm}|p{0.8cm}|p{0.7cm}|p{0.7cm}|p{0.7cm}|p{0.9cm}|p{0.77cm}|p{0.9cm}|}
			\hline
			\multicolumn{19}{|c|}{Keywords} \\
			\hline
			Author&	Poisson Arrival& MAP&	BMAP&	MMAP& Expon- ential Servi- ce Time&	PH Service& Prod- uct Form Solu- tion&	Peris- hable 
			Inven- tory& Impat- ient Customer&	Single Server&	Multi Server&	Retrial Queues&	MAM&	Conti- nuous Review&	Discr- ete QIS&	Produ- ction Inve- ntory&	Server Vaca- tion&	Multi-class cust- omers \\
			\hline
			V S Koroliuk\cite{kv} &$\checkmark$&&&&$\checkmark$&&&$\checkmark$&&&&&&&&&$\checkmark$& \\
			\hline
			K P Jose\cite{kpjp} &$\checkmark$&&&&$\checkmark$&&&&&&$\checkmark$&&$\checkmark$&&&$\checkmark$&$\checkmark$& \\
			\hline
			K Prasanna Lakshmi\cite{kpl} &$\checkmark$&&&&$\checkmark$&&&$\checkmark$&&$\checkmark$&&&$\checkmark$&&&&$\checkmark$& \\
			\hline
			A Krishna- moorthy\cite{ak} &&&$\checkmark$&&&&&&&$\checkmark$&&&&&&&$\checkmark$& \\
			\hline
			A Krishna- moorthy\cite{kavr} &$\checkmark$&&&&&&$\checkmark$&&&&&&&&&$\checkmark$&& \\
			\hline
		    K R Rejitha\cite{krr} &$\checkmark$&&&&&$\checkmark$&&&&&&$\checkmark$&$\checkmark$&&&&& \\
			\hline
			P V Laxmi\cite{lp} &&$\checkmark$&&&&&&$\checkmark$&&&&$\checkmark$&$\checkmark$&$\checkmark$&&&$\checkmark$&$\checkmark$ \\
			\hline
            Linhong Li\cite{llwz} &$\checkmark$&&&&$\checkmark$&&$\checkmark$&&&$\checkmark$&&&&&&&$\checkmark$& \\
			\hline
            Manikandan R\cite{mr} &$\checkmark$&&&&$\checkmark$&&$\checkmark$&&&&&&&&&$\checkmark$&& \\
			\hline
			K A Maqbali \cite{Maqb} &$\checkmark$&&&&&$\checkmark$&&&&$\checkmark$&&&$\checkmark$&&&&& \\
			\hline
			M P Anilkumar\cite{mpank} &&&&&&&&&&&&&$\checkmark$&&$\checkmark$&$\checkmark$&& \\
			\hline
			A Melikov\cite{maas} &$\checkmark$&&&&&&&&&&&$\checkmark$&&&&&&$\checkmark$ \\
			\hline
			A Melikov\cite{ma} &$\checkmark$&&&&&&&$\checkmark$&$\checkmark$&&&&&&&&& \\
			\hline

		\end{tabular}
	\caption{QIS with different categories}
	\end{table}
    \end{landscape}

    \begin{landscape}
	\begin{table}
		\begin{tabular}{|p{1.3cm}||p{0.9cm}|p{0.6cm}|p{0.8cm}|p{0.9cm}|p{0.4cm}|p{0.7cm}|p{0.7cm}|p{0.8cm}|p{0.75cm}|p{0.8cm}|p{0.7cm}|p{0.7cm}|p{0.9cm}|p{0.65cm}|p{0.7cm}|p{0.6cm}|p{0.8cm}|p{0.77cm}|p{0.8cm}|}
			\hline
			\multicolumn{20}{|c|}{Keywords} \\
			\hline
			Author&	Poisson Arrival& MAP&	BMAP& Expon- ential Servi- ce Time&	PH Service& Prod- uct Form Solu- tion&	Peris- hable 
			Inven- tory&	Multi- comm odity&	Post poned demand&	Impat- ient Customer&	Single Server&	Multi Server&	Retrial Queues&	MAM&	Conti- nuous Review&	Discr- ete QIS&	Produ- ction Inve- ntory&	Server Vaca- tion&	Multi-class cust- omers \\
			\hline
			A Z Melikov\cite{am} &$\checkmark$&&&$\checkmark$&&&&&&&&&&&&&&&$\checkmark$ \\
			\hline
			A Melikov\cite{measa} &$\checkmark$&&&&&&$\checkmark$&&&$\checkmark$&&&&&&&&& \\
            \hline
            N Anbazhagan \cite{anjgp} & &$\checkmark$   & &$\checkmark$ & & &$\checkmark$ &$\checkmark$ & & & & & & & & & & & \\
			\hline
			N Sangee-tha\cite{nst} &$\checkmark$&&&$\checkmark$&&&&&$\checkmark$&&&&&&$\checkmark$&&&&$\checkmark$ \\
			\hline
			N Saranya\cite{ns} &&&&&&&$\checkmark$&&$\checkmark$&&&&&&$\checkmark$&&&& \\
			\hline
			S Otten \cite{osk} &$\checkmark$&&&&&$\checkmark$&&&&&&&&&&&$\checkmark$&& \\
			\hline
			S Otten\cite{od} &$\checkmark$&&&$\checkmark$&&&&&&&$\checkmark$&&&&&&$\checkmark$&&$\checkmark$ \\
			\hline
			S Otten  \cite{os}&$\checkmark$&&&$\checkmark$&&$\checkmark$&&&&&$\checkmark$&&&&&&$\checkmark$&& \\
			\hline
			S Ozkar \cite{oz} &$\checkmark$&&&&$\checkmark$&&&$\checkmark$&&&$\checkmark$&&&&&&&&$\checkmark$ \\
			\hline
			P Beena\cite{pbj} &&$\checkmark$&&&&&&&&&&$\checkmark$&&$\checkmark$&&&$\checkmark$&$\checkmark$& \\
			\hline
			Pikkala Vijaya Laxmi\cite{pv} &$\checkmark$&&&$\checkmark$&&&&&&&&&$\checkmark$&$\checkmark$&$\checkmark$&&&$\checkmark$& \\
			\hline

		\end{tabular}
	\caption{QIS with different categories}
	\end{table}
    \end{landscape}

    \small
	\begin{landscape}
		\begin{table}
		\begin{tabular}{|p{1.8cm}||p{0.9cm}|p{0.6cm}|p{0.8cm}|p{0.8cm}|p{1cm}|p{0.4cm}|p{0.7cm}|p{0.7cm}|p{0.8cm}|p{0.8cm}|p{0.79cm}|p{0.9cm}|p{0.7cm}|p{0.7cm}|p{0.7cm}|p{0.9cm}|p{0.77cm}|p{0.9cm}|}
			\hline
			\multicolumn{19}{|c|}{Keywords} \\
			\hline
			Author&	Poisson Arrival& MAP&	BMAP&	MMAP& Expon- ential Servi- ce Time&	PH Service& Prod- uct Form Solu- tion&	Peris- hable 
			Inven- tory&	Multi- comm odity&		Single Server&	Multi Server&	Retrial Queues&	MAM&	Conti- nuous Review&	Discr- ete QIS&	Produ- ction Inve- ntory&	Server Vaca- tion&	Multi-class cust- omers \\
			\hline
			V Radhamani\cite{rvsba} &&$\checkmark$&&&$\checkmark$&&&$\checkmark$&&&&&&$\checkmark$&&&$\checkmark$& \\
			\hline
			K Rasmi\cite{rk} &&&&$\checkmark$&$\checkmark$&&&&&&$\checkmark$&&$\checkmark$&&&&&$\checkmark$ \\
			\hline
			P S Reshmi\cite{rpsj} &&&&$\checkmark$&&$\checkmark$&&$\checkmark$&&&&$\checkmark$&$\checkmark$&&&&& \\
			\hline
			R Manikandan\cite{rms} &$\checkmark$&&&&$\checkmark$&&&&&&&&&&&&$\checkmark$& \\
			\hline

            S K Samanta\cite{siv} &$\checkmark$&&&&$\checkmark$&&&&&$\checkmark$&&&&$\checkmark$&&&& \\
			\hline
   
			N Sangeetha \cite{sns} &&$\checkmark$&&&$\checkmark$&&&$\checkmark$&&&&&&$\checkmark$&&&& \\
			\hline
            J S A Jenifer\cite{sass} &&&&&&$\checkmark$&&&&$\checkmark$&&&&&&&& \\
			\hline
			Serife Ozkar\cite{so} &$\checkmark$&&&&$\checkmark$&&&&$\checkmark$&&&&$\checkmark$&&&&&$\checkmark$ \\
			\hline
			D Shajin\cite{sd} &&$\checkmark$&&&$\checkmark$&$\checkmark$&&&$\checkmark$&$\checkmark$&&&&&&&& \\
			\hline
			D Shajin\cite{sdjvvm} &&&&&$\checkmark$&&&&&$\checkmark$&&&&&&&&$\checkmark$ \\
			\hline
			D Shajin \cite{sdk} &$\checkmark$&&&&$\checkmark$&&$\checkmark$&&&$\checkmark$&&$\checkmark$&&&&&& \\
			\hline
			D Shajin \cite{skm} &$\checkmark$&&&&$\checkmark$&&$\checkmark$&&&$\checkmark$&&&&&&&& \\
			\hline
			D Shajin \cite{skd} &&$\checkmark$&&&&$\checkmark$&$\checkmark$&&&$\checkmark$&&&&&&&& \\
			\hline
			D Shajin \cite{skms} &$\checkmark$&&&&$\checkmark$&&&&&&$\checkmark$&&&&&$\checkmark$&& \\
			\hline
			S Jehoashan Kingsly\cite{sjk} &$\checkmark$&&&&$\checkmark$&&&$\checkmark$&&&$\checkmark$&&&$\checkmark$&&&&$\checkmark$ \\
			\hline
            S Otten\cite{srhk} &&$\checkmark$&&&$\checkmark$&&&&&&&&&&&&& \\
			\hline
			S R Chakravar-thy\cite{src} &&&$\checkmark$&&&$\checkmark$&&&&$\checkmark$&&&$\checkmark$&&&&& \\
			\hline
		\end{tabular}
	\caption{QIS with different categories}
	\end{table}
    \end{landscape}

    \begin{landscape}
			\begin{table}
		\begin{tabular}{|p{1.8cm}||p{1cm}|p{0.6cm}|p{1cm}|p{0.4cm}|p{0.7cm}|p{0.7cm}|p{0.8cm}|p{0.75cm}|p{0.8cm}|p{0.8cm}|p{0.79cm}|p{0.9cm}|p{0.7cm}|p{0.7cm}|p{0.7cm}|p{0.9cm}|p{0.77cm}|p{0.9cm}|}
			\hline
			\multicolumn{19}{|c|}{Keywords} \\
			\hline
			Author&	Poisson Arrival& MAP&	 Expon- ential Servi- ce Time&	PH Service& Prod- uct Form Solu- tion&	Peris- hable 
			Inven- tory&	Multi- comm odity&	Post poned demand&	Impat- ient Customer&	Single Server&	Multi Server&	Retrial Queues&	MAM&	Conti- nuous Review&	Discr- ete QIS&	Produ- ction Inve- ntory&	Server Vaca- tion&	Multi-class cust- omers \\
			\hline
			C Sughanya \cite{scs} &&$\checkmark$&&$\checkmark$&&&&&&&$\checkmark$&$\checkmark$&&&&&$\checkmark$& \\
			\hline
			C Sughanya \cite{cs} &&$\checkmark$&&$\checkmark$&&&&&&&$\checkmark$&&&$\checkmark$&&&&$\checkmark$ \\
			\hline
			C Sugapriya \cite{scnm} &&&&&&&&&&$\checkmark$&&$\checkmark$&$\checkmark$&$\checkmark$&&&$\checkmark$& \\
			\hline
			H Tayebi \cite{th} &$\checkmark$&&$\checkmark$&&&&&&&&&&&&&&& \\
			\hline
			U Uzunoglu Kocer \cite{uk} &$\checkmark$&&&&&$\checkmark$&&&&&&&&$\checkmark$&&&&$\checkmark$ \\
			\hline
			D T Varghese \cite{vd} &$\checkmark$&&$\checkmark$&&&&&&&$\checkmark$&&&$\checkmark$&&&&& \\
            \hline
			V Vinitha \cite{vaaj} &&$\checkmark$&$\checkmark$&&&&&&&&&&&&&&&$\checkmark$ \\
			\hline
			FF Wang\cite{wfb} &&$\checkmark$&&&&&&&&&$\checkmark$&$\checkmark$&&&&&&$\checkmark$ \\
			\hline
			Y Barron\cite{ybar} &&&&&&&&&&&&&&$\checkmark$&&$\checkmark$&& \\
			\hline
			Y Barron\cite{ybo} &$\checkmark$&&$\checkmark$&&&$\checkmark$&&&&&&&&$\checkmark$&&&& \\
			\hline
			D Yue \cite{yq} &$\checkmark$&&$\checkmark$&&$\checkmark$&&&&&$\checkmark$&&&&&&$\checkmark$&$\checkmark$& \\
			\hline
			D Yue \cite{ywz} &$\checkmark$&&$\checkmark$&&$\checkmark$&$\checkmark$&&&&$\checkmark$&&&&&&$\checkmark$&$\checkmark$& \\
			\hline
			D Yue \cite{yzq} &$\checkmark$&&&&$\checkmark$&&&&&$\checkmark$&&&&&&&& \\
			\hline
			D Yue \cite{yzy} &$\checkmark$&&$\checkmark$&&&&&&&&$\checkmark$&&&&&&& \\
			\hline
			Y Zhang\cite{yzdy} &$\checkmark$&&$\checkmark$&&&&&&$\checkmark$&&&&&$\checkmark$&&&$\checkmark$& \\
			\hline
			Y Zhang \cite{zy} &$\checkmark$&&$\checkmark$&&$\checkmark$&&&&&&&&&$\checkmark$&&&$\checkmark$& \\
			\hline

		\end{tabular}
  
	\caption{QIS with different categories}
	\end{table}
	\end{landscape}

 \begin{center}
    \begin{table}\small
	\vspace{1cm}
	\begin{tabular}{|p{1.1cm}||p{1cm}|p{1.1cm}|p{1cm}|p{1.1cm}|p{1.4cm}|p{1.1cm}|p{1cm}|}
		\hline
		\multicolumn{8}{|c|}{Keywords} \\
		\hline
		Author&Poisson Arrival &Expon-ential Service &Single Server &Server Vacation &Production Inventory &Continu-ous Review &Game Theoretic\\
		\hline
		Jinting Wang \cite{Wang1} &$\checkmark$&$\checkmark$&$\checkmark$&&&$\checkmark$&$\checkmark$ \\
		\hline
		B Kim \cite{Kim} &$\checkmark$&&&&$\checkmark$&&$\checkmark$ \\
		\hline
		D H Lee  \cite{Lee} &$\checkmark$&$\checkmark$&$\checkmark$&&&&$\checkmark$ \\
		\hline
		Q Li \cite{Guo} &$\checkmark$&&&$\checkmark$&&&$\checkmark$ \\
		\hline
		R Wang \cite{Wang2} &&&&&&&$\checkmark$ \\
		\hline
		 Xiaoli Cai \cite{Cai} &$\checkmark$&$\checkmark$&$\checkmark$&&&&$\checkmark$ \\
		\hline
		X Zhang \cite{Zhang} &$\checkmark$&$\checkmark$&$\checkmark$&&&&$\checkmark$ \\
		\hline
	\end{tabular}
    \caption{Game Theoretic QIS with different categories}
    \end{table}
    \end{center}

\textit{Acknowledgement.} Salini K's research is supported by Research supported by Kerala State Council for Science, Technology \& Environment (KSCSTE) (No. KSCSTE/2092/2019-FSHP- MAIN) and the research work of Manikandan Rangaswamy, is supported by SERB, India, grant EEQ/2022/000229.

\bibliographystyle{unsrt}  
%\bibliography{references}  %%% Remove comment to use the external .bib file (using bibtex).
%%% and comment out the ``thebibliography'' section.

%%% Comment out this section when you \bibliography{references} is enabled.

\end{document}